\documentclass[11pt]{article}
\usepackage{amssymb,amsthm,amsmath}
\usepackage[all]{xy}
\usepackage{pxfonts} 

\theoremstyle{plain}
    \newtheorem{theorem}{Theorem}
    \newtheorem{corol}[theorem]{Corollary}
    \newtheorem{prop}{Proposition}[section]
    \newtheorem{othertheorem}[prop]{Theorem}
    \newtheorem{lemma}[prop]{Lemma}
    \newtheorem{problem}[prop]{Problem}
\theoremstyle{definition}
    \newtheorem{defn}[prop]{Definition}
\theoremstyle{remark}
    \newtheorem{rem}[prop]{Remark}
    \newtheorem{question}[prop]{Question}

\newcommand{\C}{\mathbb{C}}
\newcommand{\R}{\mathbb{R}}
\newcommand{\Q}{\mathbb{Q}}
\newcommand{\Z}{\mathbb{Z}}
\newcommand{\N}{\mathbb{N}}
\newcommand{\T}{\mathbb{T}}
\newcommand{\I}{\mathbb{I}}
\newcommand{\G}{\mathbb{G}}
\newcommand{\D}{\mathbb{D}}
\renewcommand{\P}{\mathbb{P}}

\newcommand{\cU}{\mathcal{U}}
\newcommand{\cL}{\mathcal{L}}

\newcommand{\cR}{\mathcal{R}}

\newcommand{\cT}{\mathcal{T}}

\newcommand{\cP}{\mathcal{P}}
\newcommand{\cO}{\mathcal{O}}
\newcommand{\eps}{\varepsilon}
\renewcommand{\setminus}{\smallsetminus}
\renewcommand{\emptyset}{\varnothing}
\newcommand{\LE}{\mathit{LE}}
\newcommand{\id}{\mathit{id}}
\newcommand{\Aut}{\mathit{Aut}}
\newcommand{\Homeo}{\mathit{Homeo}}
\DeclareMathOperator{\diam}{diam}
\DeclareMathOperator{\trace}{tr}

\DeclareMathOperator{\interior}{int}

\begin{document}

\title{Dichotomies between uniform hyperbolicity and zero Lyapunov
exponents for $\mathit{SL}(2,\mathbb{R})$ cocycles}

\date{October, 2005. Revised: June, 2006.}

\author{Jairo Bochi\thanks{Partially supported by CNPq-Profix and
Franco-Brazilian cooperation program in Mathematics.} \ and Bassam  Fayad}



\maketitle

\begin{abstract}
We consider the linear cocycle $(T,A)$ induced by
a measure preserving dynamical system $T \colon X \to X$ and
a map $A \colon X \to \mathit{SL}(2,\mathbb{R})$.
We address the dependence of the
upper Lyapunov exponent of $(T,A)$ on the dynamics $T$ when
the map $A$ is kept fixed.
We introduce  explicit conditions on the cocycle that allow to perturb the dynamics,
in the weak and uniform topologies, to make the exponent drop arbitrarily close to zero.

In the weak topology we deduce that if $X$ is a compact connected manifold, then for a 
$C^r$ ($r \ge 1$) open and dense set of maps $A$, either
 $(T,A)$ is uniformly hyperbolic for every $T$, or  the  Lyapunov exponents of $(T,A)$
vanish for the generic measurable $T$.

For the continuous case, we obtain that if $X$ is of dimension greater than $2$, then for a 
$C^r$ ($r \ge 1$) generic map $A$, there is a residual set of 
volume-preserving homeomorphisms $T$ for which either $(T,A)$ is uniformly hyperbolic 
or the Lyapunov exponents of $(T,A)$ vanish.

\vskip 0.5\baselineskip

\vskip 0.5\baselineskip
\end{abstract}

\section{Introduction and statement of the results}

\emph{Throughout this paper let $\G = \mathit{SL}(2,\mathbb{R})$.}

Let $\mu$ be a finite positive measure on a measurable space $X$,
$T \colon X \to X$ be a $\mu$-preserving map,
and $A\colon X\to\G$ be a measurable map.
The pair $(T,A)$ is called a cocycle.
It induces a skew-product map
$F_{T,A} \colon X \times \R^2 \to X \times \R^2$ defined by
$F_{T,A}(x,v) = (T(x), A(x)v)$.
Denote, for $n \in \N$,
$$
A_T^n(x) = A(T^{n-1}(x)) \cdots A(T(x)) A(x),
\quad \text{so that} \quad
F^n_{T,A}(x,v) = (T^n(x), A^n_T(x)v).
$$
The group $\G$ also acts on $\P^1 = \P(\R^2)$,
so $(T,A)$ also induces a skew-product map on $X \times \P^1$.
By simplicity we use the same notations in both $\R^2$ and $\P^1$ cases.

Provided $\log\|A\| \in L^1(\mu)$,
the upper Lyapunov exponent of the cocycle $(T,A)$ at $x\in X$, given by
$$
\lambda(T,A,x) = \lim_{n \to \infty} \frac{1}{n}\log\|A_T^n(x)\|,
$$
exists for $\mu$-almost every $x \in X$.
(See e.g.~\cite{LArnold} for basic facts about Lyapunov exponents.)
We denote also
$$
\LE(A,T) = \int_X \lambda(T,A,x)\,d\mu(x).
$$

A cocycle $(T,A)$ where $A$ is essentially bounded
is called \emph{uniformly hyperbolic} if there exists,
for $\mu$-a.e.~$x\in X$, a splitting $\R^2 = E^u(x) \oplus E^s(x)$, 
which varies measurably with respect to $x$, is $F_{T,A}$-invariant,
and such that $E^u$ is uniformly expanded and $E^s$ is uniformly contracted.

Uniform hyperbolicity of $(T,A)$ is equivalent to the following:
there exists $c>0$, $\lambda>1$ such that
$\| A_T^n(x) \| > c \lambda^n$ for $\mu$-a.e.~$x$ and $n \ge 0$.
See \cite[proposition 2]{Yo}. 

\medskip

In this paper we address the question
of the dependence of $\LE(A,T)$ on the dynamics $T$,
where $A\colon  X \to \G$ is fixed.
We shall consider the following two general situations:

\begin{description}

\item[{\it Measurable situation:}]
Assume that $(X,\mu)$ is a non-atomic Lebesgue space
and $A \colon X \to \G$ is a bounded measurable map.
The dynamics $T$ varies in the space $\Aut(X,\mu)$
of the automorphisms of $(X,\mu)$ (i.e., bi-measurable $\mu$-preserving bijections).
We will always consider the space $\Aut(X,\mu)$
endowed with the \emph{weak topology}, according to which
$T_n \to T$ iff $\mu(T_n(B) \mathbin{\vartriangle} T(B)) \to 0$
for every measurable set $B \subset X$.

\item[{\it Continuous situation:}]
$X$~is a compact manifold of dimension at least~$2$, 
$\mu$~is a volume measure,
and \mbox{$A\colon X \to \G$} is continuous.
Now the dynamics $T$ varies in the space $\Homeo(X,\mu)$ of $\mu$-preserving homeomorphisms,
which we endow with the uniform ($C^0$) topology.
\end{description}

\begin{rem}\label{r.generic ergodicity}
In the respective topologies,
the generic maps $T \in\Aut(X,\mu)$ and $T \in \Homeo(X,\mu)$ are \emph{ergodic}.
These are classical theorems of Halmos~\cite{H}
and Oxtoby and Ulam~\cite{OU}, respectively.
Although we shall not use these results, we use some related ideas from \cite{AlPr}.
\end{rem}

\medskip

More concretely, we are interested in finding out
\emph{under what conditions on the maps $A$ and $T$, is it possible to find a perturbation}
(in one of the two topologies above, according to the situation considered)
\emph{$\tilde{T}$ of $T$ so that the Lyapunov exponent $\LE(A,\tilde{T})$
drops to an arbitrarily small value or even to zero.}
Moreover, we want those conditions to be \emph{checkable},
instead of appealing to Baire's theorem.
It is clear that such conditions must exclude some kind of hyperbolicity:
if for example $A$ is constant and hyperbolic, then $\LE(A,T)$ is positive and 
independent of $T$.
(We will see later that the ``kinds'' of hyperbolicity we have to exclude are 
not the same in the measurable and continuous situations.)
With this in mind, we look for the weakest possible conditions that imply
dichotomies between zero exponents and uniform hyperbolicity.

\medskip

In the measurable situation, we shall define a condition
over $A$, called \emph{richness},
that guarantees the existence of maps $T$ such that $\LE(A,T)=0$.
In fact, we will prove (theorem~\ref{t.meas}) that
\emph{if $A$ is rich then the generic
$T \in \Aut(X,\mu)$ satisfies $\LE(A,T)=0$.}
The richness condition is explicit and involves only the measure $\nu = A_* \mu$
(see definition~\ref{d.rich}).
It provides some ``abundance'' of matrices in the support of $\nu$
that makes it possible to find elliptic products, ``mix directions'',
and make the exponents vanish after a perturbation of the dynamics.

Richness involves absolute continuity.
So it turns out that
\emph{for the richness condition to be checkable, we have to ask some differentiability of $A$}.
Since we are working in the measurable category,
that restriction may be considered as a drawback.
However, it is perhaps  inevitable that some higher regularity must be asked from $A$.
We \emph{conjecture} indeed (see \S~\ref{ss.openquestions}) that
\emph{there exists a map \mbox{$A \colon X \to \G$}} (that assumes finitely many values only)
\emph{such that the integrated exponent $\LE(A,T)$
is bounded from below by a constant $\lambda_0>0$ for all $T\in \Aut(X,\mu)$,
and nevertheless $A$ assumes an elliptic value on a set of positive measure.}

\medskip

\begin{rem} \label{r.when A varies}
Our problem of studying $\LE(A,T)$ as a function of $T$
is parallel to the one addressed in \cite{B},
where the map $T$ on the base is fixed and the cocycle $A$ is perturbed.
The following results are obtained in~\cite{B}:
For any ergodic $T\in\Aut(X,\mu)$,
there is a residual set $\cR_T \subset L^\infty (X,\G)$
such that if $A\in\cR_T$ then either the cocycle $(T,A)$ is uniformly hyperbolic or $\LE(A,T)=0$.
If, in addition, $X$ is a compact Hausdorff space and $\mu$ is a Borel measure then
there exists a residual subset $\cR'_T \subset C^0 (X,\G)$ with the same properties.
(See also \cite{BV} and \cite{ArB} for extensions and related results.)
\end{rem}

\medskip

\begin{defn}\label{d.measdich}
Let us say that a bounded map $A \colon X \to \G$ \emph{satisfies the measurable dichotomy} if:
\begin{itemize}
\item either there is a constant $\lambda>1$ such that 
$\|A(x_n) \cdots A(x_1)\| > \lambda^n$ for every $x_1$, \ldots, $x_n \in X$
(in particular, for every $T\in\Aut(X,\mu)$ the cocycle $(A,T)$ is uniformly hyperbolic);

\item or the set of $T\in \Aut(X,\mu)$ for which $\LE(A,T)=0$ is residual in the weak topology.
\end{itemize}
\end{defn}

Our next goal would be to express that \emph{``most'' $A$'s satisfy the measurable dichotomy}.
We are able to prove that fact if we restrict ourselves
to $C^1$ maps $A \colon X \to \G$ on some connected manifold $X$.
In fact, we give a complete classification in that case
(see theorem~\ref{t.measclassif}) that permits us to: 1)~obtain the measurable dichotomy for an open and dense set of maps $A$ (corollary \ref{c.measdich}), 2)~describe precisely which maps do not satisfy the dichotomy and explicit all the possible behaviors of the Lyapunov exponents for these maps (proposition~\ref{p.addendum}).

In this regard, we also show that if we consider only maps $A$ that assume finitely many values,
then ``most'' of them will satisfy dichotomy (see theorem~\ref{t.discrete}).
However, that case is unrelated to richness,
and actually we lose any explicit condition for zero exponents.

\medskip

In the continuous setting, we show that under certain conditions
on the map $A$ and also on the volume-preserving homeomorphism $T$, it is
possible to perturb $T$ and make the exponent drop or
(under a stronger assumption) vanish.
We obtain those results as consequences of the afore-mentioned measurable results.
For that we use Alpern's technique \cite{Alpern Lusin} of approximating
measurable automorphisms by homeomorphisms.

At last, we obtain a dichotomy result concerning homeomorphisms.

\begin{defn}\label{d.contdich}
We say that a continuous map \emph{$A \colon X \to \G$ satisfies the continuous dichotomy}
if for the generic $T\in\Homeo(X,\mu)$:
\begin{itemize}
\item either the cocycle $(A,T)$ is uniformly hyperbolic;
\item or $\LE(A,T)=0$.
\end{itemize}
\end{defn}

We prove that if $X$ is a $C^r$ manifold ($r\ge 1$)
then the \emph{$C^r$-generic map $A$ satisfies the continuous dichotomy}.

Next we give the precise statements.

\subsection{The measurable case}

To define richness we need  first to introduce some notation.
If $\nu$ is a measure in $\G$ and $v \in \P^1$, then the push-forwards of $\nu$ by the maps
$$
M \in \G \mapsto M^{-1} \in \G \quad \text{and} \quad
M \in \G \mapsto M \cdot v \in \P^1
$$
are indicated by $\nu^{-1}$ and $\nu * v$, respectively.
If $n \in \N$, the push-forward of $\nu^n$ by the map
$$
(M_1, \ldots, M_n) \in \G^n \; \mapsto \; M_n \cdots M_1 \in \G
$$
(i.e., the $n$-th convolution power) is indicated by $\nu^{*n}$.

\begin{defn}\label{d.rich}
Let $\nu$ be a finite measure on $\G$ of bounded support.
Then $\nu$ is called \emph{$N$-rich}, for $N \in \N$,
if there exists $\kappa>0$ such that for every $v \in \P^1$
we have
$$
\nu^{*N} * v \ge \kappa m
\quad \text{and} \quad
(\nu^{*N})^{-1} * v \ge \kappa m,
$$
where $m$ denotes Lebesgue measure in $\P^1$.
The measure $\nu$ is called \emph{rich} if it is $N$-rich for some $N\in\N$.
\end{defn}

The richness property is studied in \S~\ref{ss.richness},
where the following criterium is obtained:

\begin{prop}\label{p.rich}
Let $M$ be a compact manifold 
(possibly $1$-dimensional, possibly with boundary, possibly not connected),
with a smooth volume measure $\mu$,
and let $A\colon M \to \G$ be a $C^1$ map.
Assume there are points $p_1$, \ldots, $p_k \in M$ such that
the matrix $A(p_k) \cdots A(p_1)$ is elliptic
and moreover $A$ is not locally constant at at least one of the $p_i$'s.
Then $A_* \mu$ is rich.
\end{prop}

Our main theorem is:
\begin{theorem}\label{t.meas}
Let $(X,\mu)$ be a non-atomic Lebesgue space.
Let $A \colon X \to \G$ be a bounded measurable map such that the measure $A_* \mu$ is rich.
Then there is a residual set $\cR_A \subset \Aut(X,\mu)$ such that
$\LE(A,T)=0$ for all $T \in \cR_A$.
\end{theorem}

The following is an informal outline of the proof.
Richness implies the existence of products that send
any chosen direction into any other.
Perturbing the dynamics in the weak topology we can
make most of the orbits periodic.
Also, we do that perturbation so that
a very small ``rich part'' of the space is kept invariant.
Then another well chosen perturbation makes most of the orbits
spend a few iterates in the rich region
in such a way that expanded directions are sent into contracting directions.
This forces the Lyapunov exponents to drop close to zero.

\medskip

From theorem~\ref{t.meas} and proposition~\ref{p.rich},
we obtain the following result:

\begin{theorem}\label{t.measclassif}
Let $X$ be a compact connected manifold (possibly $1$-dimensional, possibly with boundary),
and $\mu$ be a volume measure.
Let $A\colon X\to\G$ be a $C^1$ map.
Then:
\begin{itemize}
\item[(i)] either there is a \emph{closed} interval $I \subset \P^1$
such that $A(x) \cdot I \subset I$ for every $x\in X$;
\item[(ii)] or $\LE(A,T)=0$ for the generic $T\in \Aut(X,\mu)$ (w.r.t.~the weak topology).
\end{itemize}
\end{theorem}

It is easy to describe all the possible behaviors of
the Lyapunov exponent under alternative~(i) in the theorem -- see \S~\ref{ss.addendum}.
In fact, for an open and dense set of maps $A$ (i) implies uniform hyperbolicity for any dynamics,
so we obtain the result mentioned in the abstract:

\begin{corol}\label{c.measdich}
Let $1 \le r \le \infty$ and let
$X$ be a compact connected $C^r$-manifold.
Then the set of $A \in C^r(X,\G)$ that satisfy the measurable dichotomy is $C^r$-open and dense.
\end{corol}

We remark that $L^\infty$ and $C^0$ versions of corollary~\ref{c.measdich} are also true:

\begin{prop}\label{p.easy dich}
Let $(X,\mu)$ be a non-atomic Lebesgue space.
There exists a residual subset $\cR\subset L^\infty(X,\G)$
such that measurable dichotomy holds for every $A \in \cR$.
If, in addition, $X$ is compact Hausdorff and $\mu$ is a Borel measure then
there exists a residual subset $\cR'\subset C^0(X,\G)$
such that measurable dichotomy holds for every $A \in \cR'$.
\end{prop}

The proposition follows from the results of \cite{B} mentioned in
remark~\ref{r.when A varies}, see appendix~\ref{ss.easydich}
for details.

\begin{question}\label{q.connectedness}
Can theorem~\ref{t.measclassif} be extended in some form to non-connected manifolds?
Does corollary~\ref{c.measdich} remain true if $X$ is not connected?
\end{question}

\subsection{The continuous case}

Next we consider volume-preserving homeomorphisms as dynamics.
As explained before, it is useful to assume differentiability of the map $A$.

So now we let $X$ be a smooth compact connected manifold,
possibly with boundary, of dimension $d \ge 2$,
and let $\mu$ be a smooth volume measure.
Recall $\Homeo(X,\mu)$ indicates the set of $\mu$-preserving homeomorphisms,
endowed with the $C^0$ topology.
Our main result in that setting is:

\begin{theorem}\label{t.cont}
Let $T\in\Homeo(X,\mu)$ and $A \in C^1(X, \G)$.
Assume there is a $T$-periodic point $p = T^n(p)$ such that:
\begin{itemize}
\item $A_T^n(p)$ is elliptic;
\item $A$ is not locally constant at at least one of the points $p$, $Tp$, \ldots, $T^{n-1}p$.
\end{itemize}
Then for every $\eps>0$ there exists $\tilde{T} \in \Homeo (X,\mu)$ arbitrarily $C^0$-close to $T$
such that $\LE(A,\tilde{T})<\eps$.
\end{theorem}

In fact, it is easy to strengthen theorem~\ref{t.cont}
by demanding only the existence of periodic pseudo-orbits with
similar properties; see \S~\ref{ss.contdich}.

\medskip

Let us call a periodic point $p$ of period $n$ for an homeomorphism $T \colon X\to X$
\emph{persistent} if for every $\eps>0$ there is $\delta>0$
such that if $\tilde{T}$ is an homeomorphism $\delta$-$C^0$-close to $T$
then $\tilde{T}$ has a periodic point $\tilde{p}$ of period $n$ which
is $\eps$-close to $p$.
(For example, if $p$ is an isolated fixed point of
Poincar\'e--Lefschetz index different from $1$ then $p$ is persistent.)

Theorem~\ref{t.cont}, together with a semicontinuity property (proposition~\ref{p.semiT}),
implies:

\begin{corol}
Let $T\in\Homeo(X,\mu)$ and $A \in C^1(X, \G)$.
Assume that $T$ has a persistent periodic point
$p = T^n(p)$ such that:
\begin{itemize}
\item $A_T^n(p)$ is elliptic;
\item for some $i=0,1,\ldots,n-1$, the derivative $DA(T^i p)$ is non-zero.
\end{itemize}
Then there is a neighborhood $\cU \subset \Homeo(X,\mu)$ of $T$
and a residual subset $\cR \subset \cU$ such that
$\LE(A,\tilde{T}) = 0$ for all $\tilde T \in \cR$.
\end{corol}

\medskip

Recall definition~\ref{d.contdich}.
The $C^0$-generic map satisfies the continuous dichotomy -- this follows 
easily from \cite{B} and proposition~\ref{p.generic fubini}.
Here we extend this result to higher topologies:

\begin{theorem}\label{t.contdich}
Let $1 \le r \le \infty$ and let $X$ be a compact $C^r$-manifold, with
a smooth volume measure $\mu$.
Then there is a residual set $\cR \subset C^r(X,\G)$
such that if $A \in \cR$ then $A$ satisfies the continuous dichotomy.
\end{theorem}

\medskip

Our proof of the corollary~\ref{c.measdich}
gives an effective way to decide whether a given map $A \colon X \to \G$
satisfies the measurable dichotomy.
That is not so for our proof of the continuous dichotomy theorem~\ref{t.contdich}.
Two related questions are:

\begin{question}
Is there any analogue of theorem~\ref{t.measclassif} for the continuous case?
\end{question}

\begin{question}
In theorem~\ref{t.contdich}, can we take an open and dense set, instead of a residual one?
\end{question}

\subsection{The smooth case}

If both $T$ and $A$ are assumed to have higher regularity,
then dichotomy between uniform hyperbolicity and zero exponents
is no longer true -- 
either fixing $A$ and varying $T$ or fixing $T$ and varying $A$.
This is shown by the following two examples:

\begin{enumerate}

\item
Let $T \colon X \to X$ be a volume-preserving $C^2$ Anosov diffeomorphism.
Then the exponent is positive 
on an open and dense subset of $C^1(X,\G)$,
by results of Bonatti and Viana~\cite{BnV}.

\item
Consider Schr\"{o}dinger cocycles on the $d$-torus:
\begin{equation}\label{e.schrodinger}
S_{\lambda,V}(\theta) =
\begin{pmatrix} \lambda V(\theta) & - 1 \\ 1 & 0 \end{pmatrix},
\quad \theta \in \T^d, \ \lambda \in \R
\end{equation}
If $V(\theta)$ is a non-constant trigonometrical polynomial,
and the dynamics in the base is restricted to real analytic maps in
a neighborhood of the unit polydisc in $\C^d$,
then Herman~\cite{HermanMethode} showed that there exists a positive lower bound on the
exponent, provided $\lambda$ is greater than some $\lambda_0$.
\end{enumerate}

\medskip

In our setting, we can ask:

\begin{question}
Assume $A \colon X \to \G$ is a $C^1$ map that assumes both
elliptic and hyperbolic values.
(A concrete interesting example in the torus $X=\T^d$ is $A=S_{V,\lambda}$
given by \eqref{e.schrodinger} with $V(\theta) = \cos \theta_1$ and $\lambda \gg 1$.)
When is it possible to find a volume-preserving $C^1$ map
$T \colon X \to X$ such that $\LE(A,T)$ is exactly zero?
\end{question}

\medskip

We mention here two results obtained by Herman proving abundance of
zero exponents in the absence of uniform hyperbolicity for smooth cocycles
above uniquely ergodic diffeomorphisms of the circle.
(Here the exponents are computed above the unique invariant probability measure.)
The results are based on Baire category arguments and the method used
is to approximate the base dynamics by periodic maps and concentrate
the measure on orbits above which the product of matrices are elliptic.

Define
$$
F_I^\infty  = \{ f \in {\rm Diff}^\infty_+(\T^1); \; \rho(f) \in \R \setminus \Q\}
$$
where $\rho(f)$ denotes the rotation number of $f$.
We consider maps $A \in C^\infty_{\neq 0}(\T^1,\G)$, that is,
smooth maps that are not homotopic to a constant matrix.
The set $ \overline{F_I^\infty} \times  C^\infty_{\neq 0}(\T^1,\G)$
is a Baire space with the $C^\infty$-topology.
Then:
\begin{othertheorem}[\cite{HermanUnpub}]
There is a dense $G_\delta$ set $G \subset \overline{F_I^\infty} \times  C^\infty_{\neq 0}(\T^1,\G)$
of cocycles $(f,A)$ such that $f$ is uniquely ergodic and
$\LE(A, f,\mu_f)=0$ with $\mu_f$ the unique invariant probability of $f$.
\end{othertheorem}
Here the absence of uniform hyperbolicity is granted by the fact
that the cocycle is not homotopic to identity -- see \cite[proposition 4.2]{HermanMethode}.

The set of smooth maps that are homotopic to a constant matrix is denoted by $C^\infty_{0}(\T^1,\G)$.
Then:
\begin{othertheorem}[\cite{HermanUnpub}]
There exists a set $F \subset F_I^\infty \times C^\infty_{0}(\T^1,\G)$
whose $C^\infty$ closure $\overline{F}$ is $C^0$ dense in
the subset of non-uniformly hyperbolic cocycles  in  $F_I^\infty \times C^\infty_{0}(\T^1,\G)$,
and such that there is  a $C^\infty$ dense $G_\delta$ set $G \subset \overline{F}$
of cocycles $(f,A)$ such that $f$ is uniquely ergodic and
$\LE(A, f,\mu_f)=0$ with $\mu_f$ the unique invariant measure of $f$.
\end{othertheorem}

\subsection{The discrete case}

We return to the measurable case and consider this time
the situation where $A\colon X \to \G$ assumes a \emph{finite} number of values.
Such $A$ cannot satisfy the richness condition, so the previous results do not apply.
Nevertheless we can prove that measurable dichotomy holds generically.

\begin{defn}\label{d.unif hyp}
A bounded set $\Sigma \subset \G$ is called \emph{uniformly hyperbolic}
if there exists $\lambda>1$ such that
$$
\|A_n \cdots A_1\| > \lambda^n
\quad \text{for all $A_1,\ldots,A_n\in\Sigma$.}
$$
\end{defn}

(Notice that the first option in definition~\ref{d.measdich} 
just amounts to saying that $A(X)$ is a uniformly hyperbolic set.)

Given an $N$-tuple of matrices $\Sigma=(A_1,\ldots,A_N)$,
for simplicity we also write $\Sigma$ for the set $\{A_1,\ldots,A_N\}$.

\begin{theorem}\label{t.discrete}
Let $N \ge 2$ be an integer.
There exists a residual set $\cR \subset \G^N$ such that for every $\Sigma \in \cR$:
\begin{itemize}
\item either $\Sigma$ is uniformly hyperbolic;
\item or for every measurable map $A \colon X \to \Sigma$ which assumes
every value in $\Sigma$ on a set of positive measure, there is a residual set
$\cR_A \subset \Aut(X,\mu)$ such that $\LE(A,T)=0$ for every $T \in \cR_A$.
\end{itemize}
\end{theorem}

\begin{rem}
For the reader's information, we mention a characterization of 
uniform hyperbolicity obtained in \cite{ABY}:
A compact set $\Sigma\subset \G$ is uniformly hyperbolic iff 
there exists an open set $U \subset \P^1$ 
with finitely many connected components and $\overline{U} \neq \P^1$,
such that $\overline{A(U)} \subset U$ for each $A \in \Sigma$.
\end{rem}

\subsection{Structure of the paper}

In section~\ref{s.proofmain} we prove theorem~\ref{t.meas}.
In section~\ref{s.consequences} we prove proposition~\ref{p.rich} and then
theorem~\ref{t.measclassif}.
Section~\ref{s.continuous} deals with the continuous case and contains the proof
of theorem~\ref{t.cont} (that uses again theorem~\ref{t.meas} and proposition~\ref{p.rich}).
The proof of theorem~\ref{t.discrete} is given in section~\ref{s.discrete}
and it is independent of the other results.
In the appendix we present some technical results that are used throughout the paper.

\section{Proof of theorem~\ref{t.meas}}\label{s.proofmain}

In all this section, $(X,\mu)$ denotes a non-atomic Lebesgue space
(not necessarily with $\mu(X)=1$).

The following are roughly the main steps of the proof:
\begin{itemize}
\item In \S~\ref{ss.lambdak}, we show that given
$A\colon (X,\mu) \to \G$, the \emph{existence} of some dynamics $T\in\Aut(X,\mu)$
for which $\LE(A,T)$ is small depends only on the push-forward measure $A_* \mu$ in $\G$.
That is very useful, because it allows us to address the questions in Lebesgue spaces $(X,\mu)$ and
maps $A\colon (X,\mu) \to \G$ that are convenient for our constructions.

\item In \S~\ref{ss.existence},   we show that if a measure $\nu=A_* \mu$
is $1$-rich then there exists $T\in\Aut(X,\mu)$ for which $\LE(A,T)$ is small.

\item In \S~\ref{ss.measures}, we collect abstract lemmas
on perturbation of measures and maps. Then, in \S~\ref{ss.convolution}, we relate the convolution measure
$\nu^{*N}$ that appears in the definition of richness with a dynamical construction.

\item In \S~\ref{ss.endoftheproof}, we conclude the proof.
We apply the result from \S~\ref{ss.convolution} to obtain an induced cocycle
whose push-forward to $\G$ is a measure which is close to an $1$-rich one.
A specific argument of continuity is used to allow the use of the results
of \S~\ref{ss.existence} despite the fact that we are dealing with a measure that is only close to a $1$-rich one.

\end{itemize}

\subsection{Least exponent of order $k$} \label{ss.lambdak}

We begin introducing some notation:
If $A\colon (X,\mu)\to \G$ is a bounded measurable map, $k \in \N$,
and $T\in\Aut(X,\mu)$, let
$$
\Lambda_k (A,T) = \frac{1}{k} \int_X \log \|A^k_T\| \; d\mu \, .
$$
Observe that by subadditivity of $\int_X \log \|A^k_T\| d \mu$ we have $\LE(A,T) \leq \Lambda_k(A,T)$.

We define the \emph{least exponent of order $k$} of $A$
$$
\Lambda_k (A) = \inf_{T \in \Aut(X,\mu)} \Lambda_k (A,T).
$$

The continuity property of $\Lambda_k$ states as follows:
\begin{lemma}\label{l.cont}
Given $k\in \N$, $C>1$ and $\delta>0$, there exists $\eta>0$
with the following properties:
If $A$, $B\colon (X,\mu)\to \G$ are measurable maps with $\|A\|_\infty$, $\|B\|_\infty \le C$
and
$$
\|A-B\|_1 = \int_X \|A-B\| \; d\mu < \eta
$$
then
$$
\left| \Lambda_k(A) - \Lambda_k(B) \right| < \delta.
$$
\end{lemma}
\begin{proof}
Take $\eta = C^{-k+1}\delta$.
Fix any $T \in \Aut(X,\mu)$.
We can estimate pointwise
$$
\left| \log \|A^k_T \| - \log \|B^k_T \|\right| \le
\left| \|A^k_T \| - \|B^k_T \|\right| \le
\| A^k_T  - B^k_T \| \le
C^{k-1} \sum_{i=0}^{k-1} \|B \circ T^i - A \circ T^i\| \, .
$$
Dividing by $k$ and integrating over $X$ we obtain
$\left| \Lambda_k(A,T) - \Lambda_k(B,T) \right| < \delta$.
\end{proof}

\begin{rem} \label{rem.isomorphism}
If $A\colon X \to \G$ is measurable and bounded, and $S\colon  (X,\mu) \rightarrow (X',\mu')$ is an isomorphism,
it is clear that $\Lambda_k (A\circ S) = \Lambda_k(A)$, because
$\Lambda_k(A \circ S, T) = \Lambda_k (A,S \circ T \circ S^{-1})$ for every $T \in \Aut(X,\mu)$.
\end{rem}

In fact, we will prove in lemma~\ref{l.push} a stronger result:
\emph{$\Lambda_k(A)$ depends only on the push-forward measure $A_* \mu$.}
To prove that, we will need lemma~\ref{l.square} below.

In what follows $\I$ denotes the unit interval $[0,1]$
and $m$ denotes Lebesgue measure on $\I$ or (by abuse of notation)
on the square $\I^2$.

\begin{lemma}\label{l.square}
Let $A\colon \I \to \G$ be measurable and bounded.
Let $\pi \colon \I^2 \to \I$ be the projection on the first coordinate,
and consider the map
$$
A \circ \pi \colon (\I^2,m) \to \G.$$
Then $\Lambda_k (A \circ \pi) = \Lambda_k(A)$ for every $k \in \N$.
\end{lemma}

The idea of the proof is to approximate $\pi$ by something invertible and
to use remark~\ref{rem.isomorphism} and the continuity property from lemma~\ref{l.cont}.
\begin{proof}
It is clear that $\Lambda_k(A \circ \pi) \le \Lambda_k(A)$, because
$\Lambda_k(A \circ \pi, T \times \id) = \Lambda_k(A,T)$ for any $T \in \Aut(\I,m)$.
Fix $\delta>0$.
Let $T \in \Aut(\I^2,m)$ be such that $\Lambda_k(A \circ \pi,T) < \Lambda_k(A\circ \pi) + \delta$.
For $n\in\N$, define an isomorphism $P_n\colon  (\I^2,m) \to (\I,m)$
such that if $I \subset \I$ is a dyadic interval with $|I|=2^{-n}$ then
$P_n(I \times \I) = I$.
Then the functions $\pi$ and $P_n\colon \I^2 \to \I$ are uniformly $2^{-n}$-close.
This implies $L^1$-convergence:
$$
\lim_{n\to\infty} \| A \circ P_n - A \circ \pi \|_1 = 0.
$$
Indeed, given $\eps>0$ there exists, by Lusin's theorem, a compact set $K \subset \I$
such that $A|_K$ is continuous and $m(K^c)<\eps$.
If $n$ is large enough then for every $x,y \in K$ that are $2^{-n}$-close
we have $\|A(x)-A(y)\|<\eps$.
Let $G_n = \pi^{-1}(K) \cap P_n^{-1}(K)$; then $m(G_n^c)< 2\eps$.
Thus
$$
\int_{\I^2} \|A \circ P_n - A \circ \pi \| \, dm =
\int_{G_n} (\cdots) + \int_{G_n^c} (\cdots) <
\eps + 2 \eps \|A\|_\infty.
$$

By lemma~\ref{l.cont}, if $n$ is sufficiently large then
$\Lambda_k(A \circ P_n, T) < \delta + \Lambda_k(A \circ \pi)$.
Let $T' = P_n \circ T \circ P_n^{-1}$; then
$\Lambda_k(A, T') = \Lambda_k(A \circ P_n, T)$.
This shows that $\Lambda_k(A) < \delta + \Lambda_k(A\circ\pi)$.
Since $\delta>0$ is arbitrary, the lemma follows.
\end{proof}

Let us record for later use something we have proved:

\begin{lemma}\label{l.pn}
There exists a sequence of isomorphisms $P_n:(\I^2,m)\to (\I,m)$ such that
for any measurable bounded $A:\I \to \G$ we have $\| A \circ P_n - A \circ \pi \|_1 \to 0$
as $n \to \infty$.
\end{lemma} 

Now we can state and prove the:

\begin{lemma}\label{l.push}
Let $A, B\colon (X,\mu) \to \G$ be such that $A_* \mu = B_* \mu$.
Then $\Lambda_k(A) = \Lambda_k(B)$.
\end{lemma}

\begin{proof}
We can assume that $\mu(X)=1$.
By remark \ref{rem.isomorphism}, we can assume that $A$ and $B$ are defined over $(X,\mu)=(\I,m)$.
Since $A_* m = B_* m$, by lemma~\ref{l.thouvenot}
there is an automorphism $S$ of the square $\I^2$ such that
$A \circ \pi = B \circ \pi \circ S$.
In particular, $\Lambda_k(A \circ \pi) = \Lambda_k(B \circ \pi)$.
So, by lemma~\ref{l.square}, $\Lambda_k(A) = \Lambda_k(B)$.
\end{proof}

Based on lemma~\ref{l.push}, we can introduce the following notation:
If $\nu$ is a finite measure in $\G$ with bounded support, and $k\in\N$,
we write
$$
\Lambda_k (\nu) = \Lambda_k(A),
$$
where $A\colon (X,\mu) \to \G$ is some map, defined on
some non-atomic Lebesgue space, such that $A_* \mu = \nu$.
(Notice the existence of such $A$.)

The map $\Lambda_k$ has the following convexity properties:

\begin{lemma}\label{l.convex}
Let $\nu$, $\nu'$ be finite measures in $\G$, with bounded supports,
and let $t>0$.
Then
\begin{enumerate}
\item[(i)]  $\Lambda_k(t\nu) = t\Lambda_k(\nu)$;
\item[(ii)] $\Lambda_k(\nu + \nu') \le \Lambda_k(\nu) + \Lambda_k(\nu')$.
\end{enumerate}
\end{lemma}

\begin{proof}
The first part is obvious, so let us show the second one.
Let $I$, $I' \subset \R$ be disjoint intervals with lengths $\nu(\G)$, $\nu'(\G)$, respectively.
Take a map $A: I \sqcup I' \to \G$ such that $A_* m|_I = \nu$,  $A_* m|_{I'} = \nu'$
(where $m$ is Lebesgue measure).
Then (1) $\Lambda_k(\nu + \nu')$, (2) $\Lambda_k(\nu) + \Lambda_k(\nu')$,
are the infimum of $\Lambda_k(A,T)$ where $T$ runs over 
(1) all automorphisms $T \in \Aut(I \sqcup I', m)$,
(2) the automorphisms $T$ such that $T(I)=I$, $T(I')=I'$, respectively.
\end{proof}

\subsection{An existence result}\label{ss.existence}

If $\nu$ is a finite measure on $\G$ with bounded support, we write
$$
|\nu| = \nu(\G) \quad \text{and} \quad
\|\nu\|_\infty = \inf\{C>1; \; \|A\| < C \text{ for $\nu$-a.e.~$A\in\G$}\},
$$
where $\|\mathord{\cdot}\|$ is some fixed operator norm.

\begin{prop}\label{p.existence}
Let $C>1$, $\delta>0$, and $\sigma$ be a $1$-rich measure with  $\|\sigma \|_\infty \le C$.
Then there exists $k \in \N$ with the following properties:
If $\omega$ is a measure in $\G$ such that
$|\omega| \le 1$ and $\|\omega\|_\infty \le C$
then
$$
\Lambda_k (\omega + \sigma) < |\omega| \delta + |\sigma| \log C.
$$
\end{prop}

The proposition implies that \emph{given $A\colon (X,\mu) \to \G$, such that
$\nu = A_* \mu$ is $1$-rich, there exists $T\in\Aut(X,\mu)$ such that
$\LE(A,T)$ is small.}
Indeed, take $\delta$ small.
The measure $\sigma = \delta \nu$ is $1$-rich as well, so we can apply the
proposition to $\omega = (1 - \delta) \nu$.

The fact that $k$ depends uniformly on $\omega$, provided  $\|\omega\|_\infty \leq C$,
will be important in the proof of theorem \ref{t.meas}.

\medskip

For the reader's convenience, we will give an informal
sketch of the proof:
\begin{itemize}
\item We first deal with the case where $\omega$ is a Dirac measure
on some hyperbolic $H\in\G$.
We use $1$-richness of $\sigma$ and the abstract lemma~\ref{l.thouvenot}
to find products of length $2$ that send the expanding direction of $H$ exactly to
the contracting one.
Then we construct a dynamics so that orbits spend a long time in the
(hyperbolic) $\omega$-part of the space, then spend two iterates in the $\sigma$-part,
then return to the $\omega$-part.
This makes $\Lambda_k(\omega+\sigma)$ small for $k$ large enough.

\item To reduce the general case to the previous one,
we (essentially) approximate a given $\omega$ by a linear combination of
Dirac measures, and use lemmas~\ref{l.cont} and~\ref{l.convex}.
\end{itemize}

Given a matrix $H \in \G$, denote by $\rho(H) \in [1,\infty)$ its spectral radius.
(That is, $\rho(H) = \max(|\lambda_1|,|\lambda_2|)$ if $\lambda_1$, $\lambda_2$ are the eigenvalues of $H$.)

\begin{proof}[Proof of proposition~\ref{p.existence}]
\emph{Definition of $k$:}
Given $C>1$ and $\delta>0$, there is $\ell_0 \in \N$ with the following properties:
If $H, R \in \G$ are matrices such that
\begin{itemize}
\item $\|H\| \le C$, $\|R\| \le C^2$;
\item $H$ is a hyperbolic matrix and moreover $\rho(H) \ge e^{\delta/4}$;
\item $R(e^u) = e^s$, where $e^u$ and $e^s \in \P^1$ denote respectively
the expanding and contracting eigendirections of $H$;
\end{itemize}
then the matrix $R H^\ell$ is elliptic for every $\ell \ge \ell_0$.
To prove this fact, take a basis $\{v^u, v^s\}$ of unit eigenvectors of $H$. 
In this basis, $R$ becomes 
$\begin{pmatrix} 0 & c \\ b & d \end{pmatrix}$.
The angle between $v^u$ and $v^s$ cannot be too small, hence we can give bounds to $b$, $c$, $d$
depending only on $C$ and $\delta$.
Then $|\trace R H^\ell| = |d| \rho(H)^{- \ell} < 2$ for sufficiently large $\ell$.

Given the $1$-rich measure $\sigma$ with $\|\sigma\|_\infty \le C$,
let $\kappa>0$ be as in definition~\ref{d.rich}.
Fix an integer 
$$
\ell \ge \max (\ell_0, 2 / \kappa ).
$$
By proposition~\ref{p.unif spectral}, there exists $k' \in \N$ such that
if $E \in \G$ is an elliptic matrix with $\|E\| \le C^{\ell+2}$
then
$$
\frac{1}{p} \log \|E^p\| < \frac{(\ell+2) \delta}{2} \quad \forall p \ge k'.
$$
Fix $p \ge k'$ large enough so that defining $k = (\ell+2) p$, we have
$$
\frac{1}{k}\log\|H^k\| < \frac{\delta}{2} \quad
\text{for all $H\in \G$ s.t.~$\|H\| \le C$ and $\rho(H)<e^{\delta/4}$.}
$$
Let us verify that $k$ has the stated properties.

\medskip

\noindent \emph{First case:}
We will first prove the proposition in the case where
$\omega$ is a Dirac measure $\delta_H$ on some $H \in \G$, with $\|H\|\le C$.
We will exhibit a Lebesgue space $(X,\mu)$,
a map $A\colon X\to\G$ with $A_* \mu = \delta_H + \sigma$,
and a dynamics $T\in \Aut(X,\mu)$ such that
$\Lambda_k (A,T) < \delta/2 + |\sigma| \log C$.

If $\rho(H) < e^{\delta/4}$ (e.g.~$H$ is elliptic or parabolic),
we simply take $T = \id$, and the claim follows.

So we assume $H$ is a hyperbolic matrix, with $\rho(H) \ge e^{\delta/4}$.
Let $e^u$ and $e^s \in \P^1$ be its expanding and contracting eigendirections, respectively.
Since $\sigma$ is $1$-rich, we have
$$
\sigma * e^u \ge \kappa m, \quad \sigma^{-1} * e^s \ge \kappa m.
$$
There are measures $\sigma_1, \sigma_2 \le \sigma$ such that
$\kappa m = \sigma_1 * e^u$ and
$\kappa m = \sigma_2^{-1} * e^s$.\footnote{Given a homomorphism
$f\colon (X,\mu) \to (Y, \nu)$ and $\nu_1 \le \nu$,
define a measure $\mu_1$ in $X$ by $\frac{d\mu_1}{d\mu} = \frac{d\nu_1}{d\nu} \circ f$;
then $f_* \mu_1 = \nu_1$.}
Let $L_1 \subset J_1$, $L_2 \subset J_2$ be intervals
with $|L_i| = \textstyle{\frac 12} \kappa = \textstyle{\frac 12} |\sigma_i|$,
$|J_i| = \textstyle{\frac 12} |\sigma|$, $J_1 \cap J_2 = \emptyset$.

Choose two measurable maps $A_i\colon  J_i \to \G$ ($i=1,2$) such that
$(A_i)_* (m|_{L_i}) = \textstyle{\frac 12} \sigma_i$ and
$(A_i)_* (m|_{J_i \setminus L_i}) = \textstyle{\frac 12} (\sigma - \sigma_i)$.
By lemma~\ref{l.thouvenot}, there exists an isomorphism $S$
such that the following diagram commutes a.e:
$$
\xymatrix{
L_1 \times \I \ar[d]_S \ar[r]^{\pi} & L_1 \ar[r]^{A_1(\mathord{\cdot})e^u}       & \P^{1} \\
L_2 \times \I          \ar[r]_{\pi} & L_2 \ar[ru]_{A_2(\mathord{\cdot})^{-1}e^s} &
}
$$
That is, for a.e.~$z \in L_1 \times \I$,
$$
(A_2 \circ \pi) (S(z)) \cdot (A_1 \circ \pi)(z) \cdot e^u = e^s.
$$
Define a convenient Lebesgue space to work on:
$$
X = \I \sqcup (J_1 \times \I) \sqcup (J_2 \times \I),
$$
The measure $\mu$ on $X$ restricted to $\I$, resp.~$J_i \times \I$, is
one-, resp.~two-, dimensional Lebesgue measure.
The map $A\colon  X \to \G$ is defined as
$A = H$ on $\I$,
$A = A_1 \circ \pi$ on $J_1 \times \I$, and
$A = A_2 \circ \pi$ on $J_2 \times \I$.
Then $A_* \mu = \delta_H + \sigma$.
At last, we define the measure-preserving dynamical system $T\colon X \to X$.
Break $\I$ into disjoint intervals $I_1$, \ldots, $I_\ell$ of equal measure
$m(I_1) = 1/\ell$.
Since $1/\ell \le \kappa/2$,
we can take a set $Z \subset L_1 \times \I$ with $m(Z) = m(I_1)$.
Let $U_1\colon  I_\ell \to Z$ be an isomorphism.
We define $T$ as being the identity on
$$
X_\id = ((J_1 \times \I) \setminus Z) \cup ((J_2 \times \I) \setminus S(Z)),
$$
and in the rest as
$$
I_1 \to I_2 \to \cdots \to I_\ell
\xrightarrow{U_1} Z \xrightarrow{S} S(Z) \xrightarrow{U_2} I_1 \, ,
$$
where the unspecified arrows are translations and the isomorphism $U_2$ 
are chosen so that $T^{\ell+2}|_{I_1}$ is identity.

Let us estimate $\Lambda_k(A,T)$.

Note that 

If $z \in I_1$ then 
$$
A^{\ell+2}_T(z) = A(S \circ T^{\ell-1} (z)) A(T^{\ell-1} (z)) H^\ell,
$$
with $A(S \circ T^{\ell-1} z) A( T^{\ell -1} z) \cdot e^u = e^s$.
By our choice of $\ell$, $A^{\ell+2}_T(z)$ is elliptic for every $z \in I_1$.
In fact this holds for any $z \in X \setminus X_\id$.
Since $p \ge k'$ and $k=(\ell+2)p$, we have, for any $z \in X \setminus X_\id$,
$$
\frac{1}{k} \log \|A^k_T(z)\| = \frac{1}{(\ell+2)p} \log \|[A^{\ell+2}_T(z)]^p\| < \frac{\delta}{2} \, .
$$
On the other hand, 
$\frac{1}{k}\int_{X_\id} \log \|A_T^k\| d\mu \le \int_{X_\id} \log \|A\| d\mu \le |\sigma| \log C$.
This shows that 
$\Lambda_k(A,T) < \frac{\delta}{2} + |\sigma| \log C$,
as claimed.

\medskip

\emph{General case:}
Now let $\omega$ be any measure satisfying the hypotheses of the proposition.
Let $I_1 = [0,|\sigma|)$, $I_2=[|\sigma|,|\sigma|+|\omega|]$, and
$A\colon  I_1 \cup I_2 \to \G$ be such that
$$
A_*(m|_{I_1}) = \sigma
\quad\text{and}\quad
A_*(m|_{I_2}) = \omega.
$$
Let $\eta = \eta(k, C, \delta|\omega|/2)$ be given by lemma~\ref{l.cont}.

Let $B\colon  I_2 \to \G$ be a simple function such that
$\|B\|_\infty \le C$ and $\|A|_{I_2}-B\|_1 < \eta$.
Extend $B$ to $I_1 \cup I_2$ by taking $B=A$ on $I_1$.
We can write
$$
B_* m = \sigma + \sum_{i=1}^n t_i \delta_{H_i} = (1-|\omega|) \sigma + \sum_{i=1}^n t_i (\delta_{H_i} + \sigma),
$$
where $H_i \in \G$, $t_i \ge 0$ and $\sum_{i=1}^n t_i = |\omega| \le 1$.

By lemma~\ref{l.convex} and the case already considered,
\begin{align*}
\Lambda_k (B)
&\le (1-|\omega|)\Lambda_k(\sigma) + \sum t_i \Lambda_k (\delta_{H_i} + \sigma) \\
&\le (1-|\omega|)|\sigma| \log C
       + \left(\sum t_i\right)\left({\textstyle \frac {\delta}2} + |\sigma| \log C\right) \\
&=   {\textstyle \frac 12}|\omega| \delta + |\sigma| \log C.
\end{align*}
Since $\|A-B\|_1 < \eta$,  we obtain
$\Lambda_k (A) < |\omega| \delta + |\sigma| \log C$.
This proves the proposition.
\end{proof}


\subsection{Perturbing measures}\label{ss.measures}

In this subsection, we introduce a definition of closeness
in the space of measures which is suitable to our purposes,
and prove a couple of useful properties.

\begin{defn}  \label{d.closeness}
Let $\nu_1$, $\nu_2$ be measures in $\G$ with bounded support and same mass $|\nu_1|=|\nu_2|=a$.
Given $\eta>0$, we say that $\nu_1$ and $\nu_2$ are \emph{$\eta$-close}
if there exist $A_1,A_2\colon ([0,a],m) \to \G$ such that $(A_i)_* \mu = \nu_i$
and $\|A_1 - A_2\|_1 < \eta$.
\end{defn}

We can define a distance $d(\nu_1, \nu_2)$ as the infimum of the $\eta$
such that $\nu_1$ and $\nu_2$ are $\eta$-close in the sense above.
That this is indeed a metric follows from the lemma below:

\begin{lemma} \label{l.L1}
Let $A \colon  (X,\mu) \to \G$ and $\nu = A_* \mu$.
If $\tilde{\nu}$ is $\eta$-close to $\nu$ then
there exists $\tilde{A}\colon  (X,\mu) \to \G$ such that
$\|\tilde{A} - A\|_1 < \eta$ and $\tilde{\nu} = \tilde{A}_* m$.
\end{lemma}

\begin{proof}
Without loss of generality we assume $|\nu|=1$, and moreover, $(X,\mu) = (\I, m)$.
By assumption, there are $A_1,A_2\colon \I \to \G$
such that $(A_1)_* m = \nu$, $(A_2)_* m = \tilde{\nu}$,
and $\eta ' = \eta - \|A_1 - A_2\|_1 > 0$.

By lemma~\ref{l.thouvenot}, there exists $S \in \Aut(\I^2,m)$
such that and $A_1 \circ \pi \circ S = A \circ \pi$ a.e.
Let $P_n \colon  \I^2 \to \I$ be given by lemma~\ref{l.pn} and
choose $n \in \N$ so that $\|A \circ \pi - A \circ P_n\|_1 <\eta'$.
Define $\tilde{A} = A_2 \circ \pi \circ  S \circ P_n^{-1}$.
Then $\tilde{A}_* m = \tilde{\nu}$ and
\begin{multline*}
\|\tilde{A} - A\|_1 =   \|\tilde{A} \circ P_n - A \circ P_n\|_1
                    =   \|A_2  \circ \pi \circ S - A \circ P_n\|_1              \le \\
                    \le \|A_2 \circ \pi \circ S - A_1 \circ \pi \circ S\|_1 + \|A \circ \pi - A \circ P_n\|_1                    <   \|A_2 - A_1 \|_1                + \eta'
                    =   \eta.
\end{multline*}
\end{proof}

We will also need the following:

\begin{lemma}\label{l.approxcut}
Let $A \colon  (X,\mu) \to \G$, $\nu = A_* \mu$ and $\sigma \le \nu$.
Then for every $\eta>0$ there exists a measurable set $Y \subset X$
such that $A_* (\mu|_Y)$ is $\eta$-close to $\sigma$.
\end{lemma}

\begin{proof}
There is no loss of generality in assuming that $|\mu|=|\nu|=1$.
So we can also assume that $(X,\mu) = (\I, m)$.
Let $f\colon  \G \to \I$ be the Radon-Nikodym derivative $\frac{d\sigma}{d\nu}$.
Define $Y_0 = \{(x, t)\in X \times \I; \; 0 \le t \le f \circ A(x)\}$.
Let $P_n \colon  X \times \I \to X$ be as in lemma~\ref{l.pn},
with $n \in \N$ large enough so that $\|A \circ \pi - A \circ P_n\|_1 < \eta$.
Let $Y = P_n (Y_0)$.
Then $A_* (\mu|_Y) = (A \circ P_n)_* ((\mu \times m)|_{Y_0})$ is $\eta$-close
to $(A\circ \pi)_* ((\mu \times m)|_{Y_0})$.
The later measure equals $\sigma$.
Indeed,
\begin{multline*}
\big( (A\circ \pi)_* ((\mu \times m)|_{Y_0}) \big) (Z) =
(\mu \times m) \big( Y_0 \cap \pi^{-1} (A^{-1}(Z)) \big) \\ =
\int_{A^{-1}(Z)} f \circ A \, d\mu =
\int_Z f \, d\nu =
\sigma(Z),
\end{multline*}
for any measurable $Z \subset \G$.
\end{proof}


\subsection{Towers and convolutions}\label{ss.convolution}

In this subsection we show that convolution measures can be
approximately obtained in a dynamical way, in the following sense:

\begin{lemma} \label{l.tower}
Let $A\colon  \I \to \G$ be bounded and $\nu = A_* m$.
Given $N \in \N$ and $\eta>0$ there exists
$F \in \Aut(\I,m)$ and a set $Z \subset \I$ such that
$F^N = \id$,
the sets $Z$, $F(Z)$, \ldots, $F^{N-1}(Z)$ are disjoint,
and the measure $(A_F^N)_* (m|_Z)$
is $\eta$-close to $\frac{1}{N} \nu^{*N}$.
\end{lemma}

\begin{proof}
If $N=1$ there is nothing to prove; assume $N \ge 2$.
Let us first consider the case where $A$ has a special form,
namely there is $M \in \N$ such that $A$ restricted to each
interval $I_j = [\frac{j-1}{M}, \frac{j}{M})$, $j=1,\ldots,M$,
is constant, say equal to $A_j$. Then
$$
\nu = \frac{1}{M} \sum_{j=1}^M \delta_{A_j} \quad\text{and}\quad
\nu^{*N} = \frac{1}{M^N} \sum_{\mathbf{j} \in \{1,\ldots,M\}^N} \delta_{A_{j_N} \cdots A_{j_1}} \, .
$$

Break each interval $I_j$ into $N M^{N-1}$ disjoint intervals of equal length,
$I_{j,k}$, $k=1,\ldots,NM^{N-1}$.
Take a bijection
$$
\{1,\ldots,N\} \times \{1,\ldots,M\}^N \to \{1,\ldots,M\} \times \{1,\ldots, NM^{N-1}\}
$$
of the form $(i, \mathbf{j}) \mapsto (j_i, k(i,\mathbf{j}))$,
where $\mathbf{j}=(j_1,\ldots,j_N)$.

Write $J_{i, \mathbf{j}} = I_{j_i, k(i,\mathbf{j})}$;
then $\{J_{i, \mathbf{j}}\}_{i, \mathbf{j}}$ is a partition of $\I$.
Define $F\colon  \I \to \I$ by mapping each
$J_{1, \mathbf{j}}$ to $J_{2, \mathbf{j}}$,
$J_{2, \mathbf{j}}$ to $J_{3, \mathbf{j}}$, \ldots, and
$J_{N, \mathbf{j}}$ to $J_{1, \mathbf{j}}$ by translations.

Let
$$
Z = \bigsqcup_{\mathbf{j} \in \{1,\ldots,M\}^N} J_{1, \mathbf{j}} \, .
$$
Then $Z$, $F(Z)$, \ldots, $F^{N-1}(Z)$ are disjoint and
$(A^N_F)_* (m|_Z) = \frac{1}{N} \nu^{*N}$.
\smallskip

\noindent \emph{General case:}
Given any bounded $A\colon  \I \to \G$, define $C= {\|A\|}_{\infty}$ and assume that
$N \in \N$ and $\eta>0$ are arbitrarily chosen.
Let $\tilde{A}\colon \I \to \G$ be a bounded simple function which is $\eps(C,N,\eta)-L^1$-close to $A$
and such that $\tilde{A}$ has $M$ level sets, all with the same measure $1/M$;
where $M$ is some integer and $\eps(C,N,\eta)$ will be defined later.

Take $S \in \Aut(\I,m)$ that maps these level sets to intervals,
so $\hat{A} = \tilde{A} \circ S$ falls in the later case.
Accordingly there exist $\hat{F} \in \Aut(\I,m)$ and a set $\hat{Z} \subset \I$
such that $\hat{Z}$, $\hat{F}(\hat{Z})$, \ldots, $\hat{F}^{N-1}(\hat{Z})$ are disjoint and
$(\hat{A}_{\hat{F}}^N)_* (m|_{\hat{Z}}) = \frac{1}{N} (\hat{A}_* m)^{*N}$.

Let $F = S^{-1} \circ \hat{F} \circ S$ and $Z = S^{-1}(\hat{Z})$.
Then $(\tilde{A}_F^N)_* (m|_Z) = \frac{1}{N} (\tilde{A}_* m)^{*N}$.
From point $2$ of the proof of lemma \ref{l.cont}
we see that $\eps(C,N,\eta)$ can be chosen so that the $\eps(C,N,\eta)-L^1$
closeness of $\tilde{A}$ and $A$ implies that  $(A_* m)^{*N}$ and $(\tilde{A}_* m)^{*N}$
are $\eta/2$-close (in the sense of definition~\ref{d.closeness}) as well as
$(A_F^N)_* (m|_Z)$ and $(\tilde{A}_F^N)_* (m|_Z)$. This concludes the proof.
\end{proof}


\subsection{End of the proof} \label{ss.endoftheproof}

An \emph{interval permutation of rank $M$} is an automorphism
$T \in \Aut(\I,m)$ which sends each interval
of the form $\left[ \frac{j}{M}, \frac{j+1}{M} \right)$, $j=0,1,\ldots,M-1$,
onto another by an ordinary translation.
We call $T$ \emph{cyclic} if the induced permutation of
the intervals is cyclic.
We are going to use the following fact:
\begin{othertheorem}[Halmos]\label{t.halmos}
Cyclic interval permutations are dense in $\Aut(\I,m)$, in the weak topology.
\end{othertheorem}
For the proof, see Halmos~\cite[p.~65]{H},
or \cite[lemmas 6.4 and 3.2]{AlPr}.

\begin{proof}[Proof of theorem~\ref{t.meas}]
Since we are working in the measurable category we can assume that
$X$ is the unit interval $\I$ and $\mu$ is the Lebesgue measure on it.

Suppose $A\colon \I \to\G$ is such that $\nu = A_* m$ is rich.
Due to proposition~\ref{p.semiT}, we only have to show
that given $T \in  \Aut(\I,m)$ and $\delta>0$,
there exists  $\tilde{T} \in  \Aut(\I,m)$ arbitrarily close to $T$
in the weak topology such that $\LE(A,\tilde{T}) < \delta$.

So let $T$ and $\delta$, and also an arbitrary $\eps>0$.
By theorem~\ref{t.halmos},
we can assume that $T$ is a cyclic interval permutation and assume
its rank $M $ satisfies $M \geq {4 / \eps}$.
Let $N$ be such that $\nu^{*N}$ is $1$-rich.

\medskip

Before going into the details of the proof let us sketch how we will
obtain the perturbation $\tilde{T}$ that will actually satisfy $m[\tilde{T} \neq T]<\eps$.
The perturbation is done in two steps.
In the first one, we will use richness of the measure $\nu$ to produce a map
$T_1$ that is close to $T$ and that has two cyclic towers:
(1) a (big) cyclic tower of height $M$ that fills most of the space and that comes
from the original tower of $T$;
(2) a (small) cyclic tower of height $N$ such that the push forward of the measure
on its base by the product of $A$'s along its $N$ levels is a measure close to
a $1$-rich one, namely $\sigma / N$, where $\sigma = \eps' \nu^{*N}$ and $\eps'$ is small.
(Actually, there is a third invariant set that can be disregarded
because it has small measure.)

Let $W$ be the union of the basis of the big and the small tower.
Then we consider the first return on the set $W$:
we obtain a derived cocycle over $W$ with identity for dynamics and a matrix map
$\hat{A}$ such that $\hat{A}_*(m|_W)$  contains a part that is close to $\sigma /N$.

Here we pass to the second step and perturb $T_1$ to $\tilde{T}$,
keeping the two towers above $W$ invariant but modifying the first return map
to a map $S$ so that the Lyapunov exponent $\LE(\hat{A},S)$ of the derived cocycle
becomes small (this is done by taking $\tilde{T}$ equal to $T_1$ except on $T_1^{-1}W$
with nevertheless $\tilde{T} T_1^{-1}W =W$).
Since $\bigcup_{n \in \Z} \tilde{T}^n W  = \bigcup_{n \in \Z} {T}_1^n W$ has almost full measure,
the latter implies smallness of $\LE(A,\tilde{T})$.
The map $\tilde{T}$ that we obtain is close to $T_1$ since we only modify the dynamics on $T_1^{-1}W$.

To understand how the map $S$ is obtained, replace for a moment $\hat{A}_*(m|_W)$ by a map
$\tilde{A}$ such that $\tilde{A}_*(m|_W)$ contains $\sigma / N$
so that proposition \ref{p.existence} applies and identity on $W$
can be replaced by a dynamics that reduces $\Lambda_k(\tilde{A},S)$ close to zero
(for some $k$ that depends on $\sigma /N$).
Now, the fact that $k$ depends only on the $1$-rich part of the measure and a
careful choice of quantifiers allow to use the continuity of $\Lambda_k$ and
derive the same conclusion for $\hat{A}$ instead of $\tilde{A}$.
Now we give the exact proof.

\medskip

Let $C = \|A\|_\infty$ and
$$
\eps' = \min \left( \eps/(4M), \delta / (M \log C) \right) \, .
$$
Let $\sigma = \eps' \nu^{*N}$.
We will use that the measure $\sigma/N$ is $1$-rich.
Let $k = k(\sigma/N,\delta,C^M)$ be given by proposition \ref{p.existence},
and let $\eta = \eta(k,\delta,C^{\max (M,N)})$ be given by lemma~\ref{l.cont}.

Using lemma~\ref{l.tower},
we find $F \in \Aut(\I,m)$ such that $F^N = \id$,
and a set $Z \subset \I$ such that $Z$,\ldots, $F^{N-1}(Z)$ are disjoint
and $(A_F^N)_*(m|_Z)$ is $\eta$-close to $\frac{1}{N} \nu^{*N}$.
By lemma~\ref{l.approxcut}, there exists a set $Y \subset Z$ such that
$(A_F^N)_*(m|_Y)$ is $\eta$-close to $\frac{1}{N} \sigma$
(from definition \ref{d.closeness} this requires that $m(Y) = |\sigma|/N$).

Let $\cT_{F} = \bigsqcup_{i=0}^{N-1} F^i(Y)$;
this set is an $F$-tower of height $N$, and has small measure:
$$m(\cT_{F}) = N m(Y) = |\sigma| < \eps'.$$

Let $\cT_T = \bigcap_{i=0}^{M-1} T^{-i}(\cT_{F}^c)$. This set has almost full measure: $m(\cT_T^c) \le M m(\cT_{F}) < M\eps' < \eps / 4$. It is also invariant by $T$ (since $T^M = id$) and we can write it as a $T$-tower of height $M$ over $I \cap {\mathcal  T}_T$ where $I$ is any interval of the cyclic permutation $T$.

Consider a first perturbation of $T$:
$$
{T}_1(x) =
\begin{cases}
T(x) & \text{if $x \in \cT_T$,} \\
F(x) & \text{if $x \in \cT_{F}$,} \\
x    & \text{otherwise.}
\end{cases}
$$
Then $\cT_T$ and $\cT_{F}$ are two disjoint invariant towers for $T_1$ with heights $M$ and $N$ respectively, and basis $I \cap \cT_T$ and $Y$ respectively. We define $\cT =  \cT_T \cup \cT_{F}$, and $W = (I \cap \cT_T) \sqcup Y$.

The first-return map of $T_1$ to $W$ is the identity.
The return-time function is $n_{W}(x)=M$ for  $x \in I \cap \cT_T$ and   $n_{W}(x)=N$ for $x \in Y$. Hence we define on $W$ the following map:
$$
\hat{A}(x) =
\begin{cases}
A_T^M(x) &\text{if $x \in I\cap \cT_T$,} \\
A_F^N(x) &\text{if $x \in Y$.}
\end{cases}
$$

Because $(A_F^N)_*(m|_Y)$ and $\sigma/N$ are $\eta$-close,
lemma~\ref{l.L1} gives a map
$\tilde{A}\colon  Y \to \G$ such that $\tilde{A}_*(m|_Y) = \sigma/N$
and $\|\tilde{A} - A_F^N|_Y\|_1 \le \eta$. Let $\hat{B} \colon  W \to \G$ be such that
$$
\hat{B}(x) =
\begin{cases}
A_T^M(x)     &\text{if $x \in I\cap \cT_T$,} \\
\tilde{A}(x) &\text{if $x \in Y$.}
\end{cases}
$$

By proposition~\ref{p.existence}, since we took $k=k(\sigma/N,\delta,C^M)$, we get
$$
\Lambda_k(\hat{B}) =
\Lambda_k (\hat{B}_* (m|_W))
  = \Lambda_k \left((A_T^M)_*(m|_{I \cap \cT_T}) + \frac{\sigma}{N}\right)   < \frac{\delta}{M} + \frac{\eps'}{N} \; \log C^M
  < 2 \delta.
$$
Since $\| \hat{A} - \hat{B} \|_1 \le 2\eta$, with $\eta =  \eta(k,\delta,C^{\max (M,N)})$,
we conclude by lemma~\ref{l.cont},
$$
\Lambda_k(\hat{A}) < 3\delta.
$$

This means that there exists an automorphism $S\colon W \to W$ such that $\Lambda_k(\hat{A}, S) < 3\delta$, and consequently $\LE(\hat{A},S) < 3 \delta$.

Finally, we define $\tilde{T}$ on $\I$:
$$
\tilde{T}(x) =
\begin{cases}
S(T_1(x))   & \text{if $x \in {T}_1^{-1}(W)$,} \\
T_1(x)      & \text{otherwise.}
\end{cases}
$$

The set $\cT=\cT_T \cup \cT_{F}$ is still invariant by $\tilde{T}$, the return time to the set $W$ is still the function $n_W$  as for $T_1$  and the products of matrices above $W$ before the first return are still given by $\hat{A}$. But the first return map to $W$ by $\tilde{T}$ is now $S$. Hence, by
proposition~\ref{p.derived}, we have $\LE(A|_\cT, \tilde{T}|_\cT)=  \LE(\hat{A}, S)$.

Recall that
$m(\cT^c) \le m(\cT_T^c) \le M \eps'$ and that we took $\eps' \leq  \delta / (M \log C)$,
therefore
$$
\LE(A,\tilde{T}) = \LE(\hat{A}, S) + \LE(A|_{\cT^c}, \tilde{T}|_{\cT^c}) < 4 \delta.
$$
We have
$$
m[\tilde{T} \neq T] \le m[\tilde{T} \neq {T}_1] + m[{T}_1 \neq T]
                    \le m(W) + m(\cT_T^c)
                    \le \left(\frac{1}{M} + \frac{\eps'}{N} \right) + \frac{\eps}{4}
                    <   \eps,
$$
as required.
\end{proof}

\section{Applying theorem \ref{t.meas}}\label{s.consequences}

\subsection{Proof of the richness criterium}\label{ss.richness}

Before deriving any consequences of theorem~\ref{t.meas},
we have to prove proposition~\ref{p.rich}.

Let us (temporarily) call a measure $\nu$ \emph{$N$-f-rich}
if there is $\kappa>0$ such that $\nu^{*N} * v \ge \kappa m$
for every $v \in \P^1$.
(The difference, compared to definition~\ref{d.rich}, is that we only consider
\underline{f}orward iterates.) 
We call $\nu$ \emph{f-rich} if it is $N$-f-rich for some $N$.

In this paragraph, $M$ denotes a compact manifold, $\interior M = M \setminus \partial M$,
and $\mu$ is a smooth volume measure on $M$.

The following lemma essentially reduces the problem of proving
richness to a one-dimensional case:

\begin{lemma}\label{l.1d reduction}
Let $A \colon M \to \G$ be a $C^1$ map.
Let $\xi \colon [-2,2] \to \interior M$ be a $C^1$ embedded arc.
Assume that
for any $v \in \P^1$,
$$
\big\{ A(\xi(t)) \cdot v ; \; t \in [-1,1] \big\} = \P^1
$$
and
$$
\frac{\partial}{\partial t} \left[ A(\xi(t))\cdot v \right] \neq 0 \quad \forall t \in [-2,2].
$$
Then $A_* \mu$ is $1$-$f$-rich.
\end{lemma}

\begin{proof}
Consider a (normal) tubular neighborhood (map) of the arc $\xi$:
$$
\Xi \colon [-2,2] \times \D^{d-1} \to M \,
$$
such that $\Xi(t,0) = \xi(t)$.
($\D^{d-1}$ denotes the open unit ball in $\R^{d-1}$.)
For $0 < \delta \le 1$ and $j = 1,2$, let
$U^j_\delta = \Xi \left([-j,j] \times \delta\D^{d-1} \right)$.
We are going to show that for sufficiently small $\delta > 0$
there exists $\kappa>0$ such that
\begin{equation}\label{e.claim}
\frac{d \left( A_* (\mu|_{U^2_\delta})* v \right)}{dm_{\P^1}}(w) \ge \kappa
\quad \forall v,w \in \P^1.
\end{equation}
We push the euclidian metric in $\R^d$ forward by
the $C^1$-diffeo $\Xi$, and without loss,
assume that $\mu|_{U^2_1}$ is the Riemannian volume induced
by that metric.
In the new metric, we have $\|\xi'(t)\|=1$.

For convenience of notation, let $v, w \in \P^1$ be fixed.
Let $F \colon M \to \P^1$ be given by $F(x) = A(x) \cdot v$.
Let $\{t_1 < \ldots < t_k\} = \{t\in [-1,1]; \; F(\xi(t)) = w\}$.
Let $D_i^\delta \subset U_\delta^2$ be the connected component of
$F^{-1}(w) \cap U_\delta^2$ that contains $\xi(t_i)$.
From the assumption, $F$ has no critical points on a neighborhood of
the are $\xi([-2,2])$.
So, for small $\delta$, $D_i^\delta$ is a $(d-1)$-submanifold.
We have
$$
\frac{d \left( F_* (m|_{U_\delta^2}) \right) }{d m_{\P^1}} (w) =
\int_{F^{-1}(w) \cap U^2_\delta} \frac{d\sigma (x)}{\| \nabla F (x)\|} \ge
\sum_{i=1}^k \int_{D_i^\delta} \frac{d\sigma (x)}{\| \nabla F (x)\|} \, ,
$$
where $\sigma$ is Riemannian $(d-1)$-dimensional volume.
Using that $\sigma (D_i^\delta) \ge \sigma (\delta \D^{d-1})$
and that $\|\nabla F\|$ is bounded, we have
$$
\int_{D_i^\delta} \frac{d\sigma (x)}{\| \nabla F (x)\|}
= \frac{\sigma(D^\delta_i)}{\|(\nabla F)(\xi(t_i)) \|} + \cO(\delta^d)
\ge C \delta^{d-1}+ \cO(\delta^d) \, ,
$$
for some $C>0$.
All bounds are uniform in $v$ and $w$, hence we can find a
small $\delta$ so that \eqref{e.claim} holds.
\end{proof}

In what follows, we denote:
$$
R_\theta =\begin{pmatrix}
\cos \theta & - \sin \theta \\
\sin \theta & \cos \theta
\end{pmatrix}.
$$

\begin{rem}\label{r.conjugacy}
For each $T\in(-2,2)$, the set of $A \in \G$ such that $\trace A = T$
consists on two $\G$-conjugacy classes, namely
$\{L R_{\theta} L^{-1}; \; L \in \G \}$ and
$\{L R_{-\theta} L^{-1}; \; L \in \G \}$,
for a certain $\theta \in (0,\pi)$.
(Note also that the set consists on a single
$\mathit{GL}(2,\R)$-conjugacy class.)
If $A$ and $B$ are two elliptic matrices in the same $\G$-conjugacy class
then there exists an unique $L$ of the form
\begin{equation}\label{e.conjugacy}
L =
\begin{pmatrix}
a & b \\ 0 & a^{-1}
\end{pmatrix},
\quad \text{with $a>0$,}
\end{equation}
such that $B = LAL^{-1}$.
\end{rem}

Next we prove the following special case of proposition~\ref{p.rich}:

\begin{lemma}\label{l.specialk}
Let $A\colon M \to \G$ be a $C^1$ map.
Assume there exists a point $p \in M$ such that
$A(p)$ is elliptic and $\trace A$ is not locally constant at $p$.
Then $A_* \mu$ is f-rich.
\end{lemma}

\begin{proof}
We can assume that $p \in \interior M$ and that $(\nabla \trace A)(p) \neq 0$.
Let $\xi \colon [-2,2] \to \interior M$ be an embedded $C^1$ path such that
$\xi(0) = p$ and
$$
\lvert \trace A(t) \rvert < 2 \quad\text{and}\quad
\frac{d}{dt}\trace A(\xi(t)) \neq 0 \quad \forall t.
$$
Let $\theta(t) = \arccos \left( {\textstyle \frac 12} \trace A(\xi(t)) \right)$.
Then $\theta(t)$ is a $C^1$ function with $\theta'(t) \neq 0$ for all $t$.
For each $t$ there exists $L(t) \in \G$ such that
$A(\xi(t))$ equals either
$L(t) R_{\theta(t)} L(t)^{-1}$ or $L(t) R_{-\theta(t)} L(t)^{-1}$
(see remark~\ref{r.conjugacy}).
By continuity, the same alternative, say the first,
occurs for all $t\in [-2,2]$.
Besides, we can choose $L(t)$ so that
$L \colon [-2,2] \to \G$ is a $C^1$ map.

For $N \in \N$, we have
$
A(\xi(t))^N = L(t) R_{N\theta(t)} L(t)^{-1} \, .
$
It is easy to see that if $N$ is large enough then
for any $v \in \P^1$,
$$
\big\{ A(\xi(t))^N \cdot v ; \; t \in [-1,1] \big\} = \P^1
\quad \text{and} \quad
\frac{\partial}{\partial t} \left( A(\xi(t))^N \cdot v \right) \neq 0 \ \forall t \in [-2,2].
$$

Next define
\begin{equation}\label{e.power}
\begin{cases}
\hat{M} = M^N = M \times \cdots \times M, \\
\hat{\mu} = \mu^N = \mu \times \cdots \times \mu, \\
\hat{A} \colon \hat{M} \to \G \text{ by } \hat{A}(x_1, \ldots, x_N) = A(x_1) \cdots A(x_N),
\end{cases}
\end{equation}
and $\hat{\xi} \colon [-2,2] \to \interior \hat{M}$ by $\hat{\xi}(t) = (\xi(t), \ldots, \xi(t))$.
Applying lemma~\ref{l.1d reduction} to these data,
we obtain that $\hat{A}_* \hat{\mu} = (A_* \mu)^{*N}$ is $1$-f-rich,
that is, $A_* \mu$ is $N$-f-rich.
\end{proof}

For the last part, we will need the following property about traces:

\begin{lemma}\label{l.trace}
Let $A$, $B\in \G$ be elliptic matrices
that are conjugate via a matrix in~$\G$.
Then
$$
\trace AB \le \trace A^2,
$$
with equality if and only if $A=B$.
\end{lemma}
The reader may check the lemma would be false
had we assumed only that the elliptic matrices have the same trace.

\begin{proof}
Write $B = L A L^{-1}$, with $L \in \G$.
We can assume that $A=R_\theta$ and
$L$ is given by~\eqref{e.conjugacy}.
Of course, $\sin \theta \neq 0$.
Direct calculation gives
$$
\trace AB =   2 - (2 + a^2 + a^{-2} + b^2) \sin^2\theta
          \le 2 - 4 \sin^2 \theta
          =   \trace A^2,
$$
with equality if and only if $a=1$ and $b=0$.
\end{proof}

Now we conclude the:

\begin{proof}[Proof of proposition~\ref{p.rich}]
It clearly suffices to show that $A_* \mu$ is f-rich.

Let $\hat{M}$, $\hat{\mu}$, and $\hat{A}$ be as in~\eqref{e.power}
with $k$ in the place of $N$.
Let also $\hat{p} = (p_1, \ldots, p_k)$.
By assumption, $\hat{A}(\hat{p})$ is elliptic
and $\hat{A}$ is not locally constant at $\hat{p}\in\hat{M}$.

If $\trace \hat{A}$ is not locally constant at $\hat{p}$ in $\hat{M}$ then,
by lemma~\ref{l.specialk},
$\hat{A}_*\hat{\mu} = (A_* \mu)^{*k}$ is f-rich, and therefore so is $A_* \mu$.

Assume then that $\trace \hat{A}$ is constant at a neighborhood of $\hat{p}$ in $\hat{M}$.
By continuity, all $\hat{A}(\hat{x})$, with $\hat{x}$ close to $\hat{p}$,
belong to a same $\G$-conjugacy class (see remark~\ref{r.conjugacy}).
Consider
$$
\check{M} = \hat{M} \times \hat {M}, \quad
\check{\mu} = \hat{\mu} \times \hat{\mu}, \quad
\check{p} = (\hat{p},\hat{p}), \quad \text{and} \quad
\check{A}(\hat{x}_1, \hat{x}_2) = \hat{A}(\hat{x}_1) \hat{A}(\hat{x}_2).
$$
Then $\trace \check{A}$ is not locally constant at $\check{p} \in \check{M}$.
(Otherwise, by lemma~\ref{l.trace}, $\hat{A}$ would be locally constant
at $\hat{p}$.)
Applying lemma~\ref{l.specialk} to $\check{A}$ we get that
$A_* \mu$ is $2k$-f-rich.
\end{proof}


\subsection{Proof of theorem~\ref{t.measclassif}}

We will need the following:

\begin{lemma}\label{l.interval}
Let $\Sigma \subset \G$ be a compact connected set.
Assume that there is no closed interval $I \subsetneqq \P^1$
such that $A \cdot I \subset I$ for every $A\in\Sigma$.
Then there are $A_1,\ldots,A_n\in \Sigma$ such that
$A_1 \cdots A_n$ is elliptic.
\end{lemma}

\begin{proof}
We claim that there is $n_0 \in \N$ such that for all $v$, $w\in\P^1$,
there exist $A_1,\ldots, A_{n_0} \in \Sigma$ such that $A_1 \cdots A_{n_0} v = w$.

Indeed, fix any matrix $A_0 \in \Sigma$, and let
$v_0 \in \P^1$ be such that $A_0 v_0 = v_0$.
Let $I_n \subset \P^1$ be the set of directions
$A_1 \cdots A_n(v_0)$, with $A_i \in \Sigma$.
Since $\Sigma$ is connected, each $I_n$ is an interval or the circle.
Also, $I_n \subset I_{n+1}$, because $v_0$ is invariant by a matrix in $\Sigma$.
Let us see that $I_{n_1} = \P^1$ for some $n_1$.
Assume the contrary, and let $I = \bigcup_n I_n$.
We have $A(\bar I) \subset \bar I$ for all $A$.
Since we are assuming $\Sigma$ has no invariant interval,
we must have $\bar I = \P^1$.
Therefore $I = \P^1 \setminus \{z\}$ for some $z$.
By the same assumption, there must be $A \in \Sigma$
such that $A(z)\neq z$.
Then $A^{-1}(z)\in I$ and so there must exist $A_1$, \ldots, $A_n$ such
that $A_1 \cdots A_n(v_0) = A^{-1}(z)$.
But this implies $z \in I$, a contradiction.

We have shown that there is $n_1$ such that for any $w$ there is
a product of length $n_1$ which sends $v_0$ to $w$.
The same reasoning applied to the set $\Sigma^{-1}$
(which does not have an invariant interval as well)
gives that there is $n_2$ such that for any $v$ there is a product of length $n_2$
sending $v$ to $v_0$.
Let $n_0 = n_1 + n_2$.
The claim is proved.

\smallskip

Next we construct an elliptic product.
Fix any $A\in \Sigma$, $A\neq \id$.
If $A$ is elliptic, we are done.

If $A$ is hyperbolic, let $e^u$ and $e^s$ be
its respectively expanding and contracting eigenvectors.
Let $B$ be a product such that $B(e^u) \in \R e^s$.
Then a calculation shows that $\trace A^n B \to 0$ as $n \to \infty$,
so there exists an elliptic product.

If $A$ is parabolic then, relative to some basis $\{e_1,e_2\}$,
$$
A = \pm \begin{pmatrix} 1 & \beta \\ 0 & 1 \end{pmatrix},
\quad \text{with $\beta \neq 0$.}
$$
Let $B$ be a product such that $B e_1 \in \R e_2$.
Write
$$
B = \begin{pmatrix} 0 & b \\ c & d \end{pmatrix},
\quad \text{with $c \neq 0$.}
$$
Then $\lvert \trace A^n B \rvert = \lvert \beta c n + d \rvert \to \infty$ as $n \to \infty$.
This shows that $\Sigma$ has a hyperbolic product.
Then we can repeat the previous reasoning and find an elliptic product.
\end{proof}

\begin{rem}
We ignore how to extend lemma~\ref{l.interval} to non-connected sets~$\Sigma$,
and that is why we are unable to answer question~\ref{q.connectedness}.
\end{rem}

We are ready now to give the:
\begin{proof}[Proof of theorem~\ref{t.measclassif}]
If the function $A$ is constant then the first case holds if $A$ is hyperbolic or parabolic,
and the second case holds if $A$ is elliptic.
So we can assume $A$ is not constant.

Assume the first case does not hold.
Applying lemma~\ref{l.interval} to $\Sigma = \{A(x) ; \; x\in X\}$,
we conclude that there is $n \in \N$ such that the function
$$
(x_1, \ldots, x_n) \in X^n \mapsto A(x_n) \cdots A(x_1)
$$
assumes an elliptic value.
This function is not constant and $X^n$ is connected,
so proposition~\ref{p.rich} applies and the exponent vanishes generically by theorem~\ref{t.meas}.
\end{proof}

\subsection{Addendum to theorem~\ref{t.measclassif} and proof of corollary~\ref{c.measdich}}
\label{ss.addendum}

We study how the Lyapunov exponent depends on $T$
if the first alternative in theorem~\ref{t.measclassif} holds.
There are two initial possibilities:
\begin{enumerate}
\item[(i.1)] There is an invariant point
(i.e. $\exists$ $v_0 \in \P^1$ s.t.~$A(x) \cdot v_0=v_0$ $\forall x\in X$).
\item[(i.2)] There is no invariant point.
\end{enumerate}
In case (i.2), consider an invariant closed interval $I$.
We further subdivide case (i.2) into cases
\emph{(i.2.1)}, \emph{(i.2.2)}, and \emph{(i.2.3)}
according to whether $A(x) \cdot I \subset I^\circ$
for \emph{none}, \emph{some but not all}, or \emph{all} $x\in X$, respectively.
(Notice that the case in which $A$ falls may depend on the choice of $I$.)

\begin{prop}[Addendum to theorem~\ref{t.measclassif}]\label{p.addendum}
Each case has the corresponding consequence as below:
\newlength{\muitotecnico} \setlength{\muitotecnico}{\textwidth}
\newlength{\maistecnico}  \settowidth{\maistecnico}{\textit{(i.2.3)}\;$\Rightarrow$\;}
\addtolength{\muitotecnico}{-\maistecnico}

\noindent
\begin{tabular}{@{}r@{\;$\Rightarrow$\;}p{\muitotecnico}@{}} 
(i.1)&There is $\lambda_0 \ge 0$ such that $\LE(A,T) = \lambda_0$ for all $T\in \Aut(X,\mu)$. \\
(i.2.1)&$\LE(A,T) > 0$ for all ergodic $T\in \Aut(X,\mu)$. \\
(i.2.2)&There is $\lambda_0 > 0$ such that $\LE(A,T) \ge \lambda_0$ for all $T\in \Aut(X,\mu)$. \\
(i.2.3)&The set $A(X) \subset \G$ is uniformly hyperbolic (see definition~\ref{d.unif hyp}),
so $(A,T)$ is uniformly hyperbolic for all $T\in \Aut(X,\mu)$.
\end{tabular}
\end{prop}

\begin{proof}
In case (i.1) we have
$$
\lambda_0 = \left| \int \log \frac{\|A(x) v_0\|}{\| v_0 \|} \, d\mu(x) \right| \, .
$$

Next recall some facts about the Hilbert metric.
That is a Riemannian metric $d$ on $I^\circ$ with the property that if
$B \in \G$ satisfies $B(I) \subset I$ then there exists $\tau_B  \ge 1$ such that
$$
d(B(x), B(y)) \le \tau_B^{-1} d(x,y) \quad \forall x, y \in I^\circ \, .
$$
Besides, if $B(I) \subset I^\circ$ then $\tau_B > 1$.
Using these facts, it is not hard to prove the remaining assertions in the proposition.
\end{proof}

\begin{proof}[Proof of corollary~\ref{c.measdich}]
For an open and dense set of $A \in C^r(X,\G)$, (i) implies (i.2.3).
\end{proof}

\section{More consequences: the continuous case}\label{s.continuous}

In this section $X$ will denote a $C^1$-smooth compact connected manifold,
possibly with boundary, of dimension $d \ge 2$,
and $\mu$ will denote a smooth volume measure.

\subsection{Notations and tools}

Here we collect some results from the book~\cite{AlPr}
that we will use in the proof of theorem~\ref{t.cont}.

Let $d$ be a metric in the manifold $X$.
The uniform topology on $\Homeo(X,\mu)$ is determined by the metric
$d(T, \tilde{T}) = \sup_{x \in X} d( T(x), \tilde{T}(x))$.
If $T$ or $\tilde{T}$ (or both) are not in $\Homeo(X,\mu)$ but in $\Aut(X,\mu)$
then the distance above should be considered with essential supremum instead of $\sup$.

In the proof of theorem~\ref{t.cont}, we will make
a non-continuous perturbation of the given homeomorphism, and
then perturb again to get a homeomorphism.
For that last step we will need
Alpern's \cite{Alpern Lusin} volume-preserving version of Lusin theorem:
\begin{othertheorem}[Theorem 10.2 from \cite{AlPr}] \label{t.lusin}
Let $T \in \Homeo(X,\mu)$ and $\eps>0$.
Then there exists $\delta>0$ with the following properties:
Given $S \in \Aut(X,\mu)$,  with $d(S, T) < \delta$
and a weak neighborhood $\mathcal{W}$ of $S$,
there exists $\tilde{S} \in \Homeo(X, \mu)$,
which equals $T$ in the boundary of $X$, and satisfies
$\tilde{S} \in \mathcal{W}$ and $d(\tilde{S}, T) < \eps$.
\end{othertheorem}

The following result permits us to carry some
geometry from the cube $(\I^d,m)$ to the manifold $(X,\mu)$.
(In fact, the proof of theorem~\ref{t.cont} would become slightly easier if $(X,\mu)=(\I^d,m)$.)

\begin{othertheorem}[Theorem 9.6 from \cite{AlPr}] \label{t.brown}
There exists a continuous map $\Phi\colon \I^d \to X$, called a \emph{Brown map},
 such that:
\begin{enumerate}
\item $\Phi$ is onto;
\item $\Phi|_{\interior \I^d}$ is a homeomorphism of the interior of $\I^d$ onto its image;
\item $\Phi(\partial \I^d)$ is a closed nowhere dense set, disjoint from $\Phi(\interior \I^d)$;
\item $\Phi_* m = \mu$.
\end{enumerate}
\end{othertheorem}

For $m \in \N$, consider the partition mod $0$ of the cube $\I^d$
into $2^{dm}$ cubes of size $2^{-m}$.
The images of those cubes by the Brown map $\Phi$ form a mod $0$ partition of $X$.
Let us indicate this partition by $\cP_m$ and call its elements \emph{cubes} as well.
An automorphism $S \in \Aut(X,\mu)$ such that
the image of a cube in $\cP_m$ is mod $0$ a cube in $\cP_m$
will be called a \emph{generalized cube exchange map}.
If additionally the map $\Phi^{-1}S\Phi$ sends cubes into cubes by translations of $\R^d$,
we call $S$ a \emph{cube exchange map}.

A generalized cube exchange map induces a permutation of the set $\cP_m$ of cubes.
We express the permutation as a product of disjoint cycles;
corresponding to each cycle there is an $S$-invariant subset of $M$,
that we call a \emph{cyclic tower}.

An important perturbation result due to Lax~\cite{L}
is the following:
\begin{othertheorem}[Theorem 3.1 from \cite{AlPr}] \label{t.lax}
Let $T \in \Homeo(\I^d, m)$.
Then for any $\delta>0$ there is a cube exchange $P$ such that $d(P,T)<\delta$.
\end{othertheorem}

In fact, it was shown by Alpern~\cite{Alpern combinatorial}
that we can take $P$ in Lax theorem with a single cycle tower.
To show that he used the lemma below
(which we will also employ with a slightly different purpose):

\begin{lemma}[Lemma 3.2 from \cite{AlPr}]\label{l.combinatorial}
Given any permutation $\sigma$ of $\mathcal{N} = \{1,\ldots,N\}$,
there is a cyclic permutation $\tilde{\sigma}$ of $\mathcal{N}$
such that $|\tilde{\sigma}(j) - \sigma(j)| \le 2$ for all $j\in \mathcal{N}$.
\end{lemma}

\subsection{Proof of theorem~\ref{t.cont}}

The proof has three steps:
\begin{enumerate}
\item
We take a fine partition $\cP_m$ of $X$ and
approximate the given $T$ by a generalized cube exchange map $S_4$.
This approximation is taken with the following additional properties:
\begin{itemize}
\item $S_4$ equals $T$ in the (periodic) $T$-orbit of $p$ (above which there is the elliptic product)
and is $C^1$ in a neighborhood of it;
\item $S_4$ has two cyclic cube towers, the smaller of them
consisting of the $n$ cubes that intersect the orbit of $p$.
\end{itemize}

\item
Using theorem~\ref{t.meas}, we change $S_4$ in a set of small diameter
to make the Lyapunov exponent vanish.
The new map $S_5$ is still uniformly close to $T$.

\item
Theorem~\ref{t.lusin} provides a homeomorphism $\tilde{T}$ weakly close to $S_4$,
which by semicontinuity will have small exponent.
\end{enumerate}

We precise now these three steps.
Let $T$, $A$, and $p$ be as in the statement, and let $\eps>0$.
Let $\delta = \delta(T,\eps)$ be given by theorem~\ref{t.lusin}.
Let $\cO(p) = \{p, Tp, \ldots, T^{n-1}p\}$.
Without loss of generality, we assume that the minimum distance between different points
in $\cO(p)$ is greater than~$\delta$.

Let $\Phi$ be the Brown map given by theorem~\ref{t.brown}.
We can assume that $\cO(p) \cap \Phi(\partial \I^d) = \emptyset$
and, moreover, that $\cO(p)$ does not intersect the cube boundaries,
for any of the partitions $\cP_m$.

\paragraph{Step 1.}
Lax theorem~\ref{t.lax} provides a cube exchange $S_1$ of rank $m$
such that $d(S_1, T) < \delta/10$.
We assume that the rank is high enough so the diameter of the cubes (in $X$)
is less than $\delta/10$.
The partition into cubes is fixed from now on.

Let $C_i$ be the cube that contains $T^i p$, for $i\in \Z_n$.
Let $H$ be the cube exchange that permutes each $S_1(C_i)$ with $C_{i+1}$,
and fixes the other cubes.
Let $S_2 = H \circ S_1$; then $S_2$ has a cyclic tower which contains $\cO(p)$
and $d(S_2, S_1) < 3 \delta/10$.

Next number all cubes from $\cP_M$ in such a way that consecutive cubes share a common face.
Then delete the $n$ cubes $C_i$ and monotonically renumber the remaining cubes,
say as $C'_j$.
Since distinct $C_i$ cubes do not share a common face,
we have that if $|j-k|\le 2$ then the diameter of $C'_j \cup C'_k$
is less than $3\delta/10$.
Applying the combinatorial lemma~\ref{l.combinatorial},
we find a cube exchange $S_3$ which is $3\delta/10$-close to $S_2$,
cyclically permutes the $C_j'$ cubes
and still satisfies $S_2(C_i)=C_{i+1}$.

Take open neighborhoods $U_i \ni T^i p$ with $\overline{U_i} \subset \interior C_i$
such that there are volume-preserving $C^1$-diffeomorphisms $f_i \colon U_i \to U_{i+1}$
satisfying $f_i(T^i p) = T^{i+1} p$.
Since $C_i \setminus U_i$ and $C_{i+1} \setminus U_{i+1}$ are Lebesgue
spaces with the same measure, we can extend each $f_i$ to
a volume-preserving map $C_i \to C_{i+1}$.
Changing $S_3$ inside the $C_i$ cubes according to those maps,\footnote{This
last step would be simpler if the Brown map $\Phi$ happened to be $C^1$.}
we obtain $S_4 \in Aut(X,\mu)$ with $d(S_4, T) < 8\delta/10$.

\paragraph{Step 2.}
The generalized cube exchange $S_4$ obtained above has two cyclic towers that cover all $X$.
We select the  cubes $C_i \ni p$ and $C_1'$ as bases of those towers.
We can assume that $C_1$ and $C_1'$ share a common face.
Let $W = C_0 \cup C_0'$; then $\diam W \le 2\delta/10$.
Consider the cocycle
$(U, \hat{A}) = \left((S_4)_W, A_{S_4,W}\right)$
induced by $S_4$ on $W$.

The measure $\hat{A}_* (\mu|_W) \ge (A_{S_4})^n_* (\mu|_{C_1})$
is rich, by proposition~\ref{p.rich}.
So theorem~\ref{t.meas} yields a measurable dynamics
$\tilde{U} \colon W \to W$ such that $\LE(\hat{A}, \tilde{U}) = 0$.
Let $S_5 \in \Aut (X,\mu)$ be given by
$$
S_5(x) =
\begin{cases}
\tilde{U}(S_4(x)) & \text{if $x \in S_4^{-1}(W)$,} \\
S_4(x)            & \text{otherwise.}
\end{cases}
$$
By proposition~\ref{p.derived}, we have $\LE(A, S_5) = 0$.
Moreover, $d(S_5, T) < \delta$.

\paragraph{Step 3.}
By semicontinuity of the Lyapunov exponent (proposition~\ref{p.semiT}),
there is a weak neighborhood $\mathcal{W} \subset \Aut (X,\mu)$ of $S_5$ such that
$LE(A, \cdot) < \eps$ on $\mathcal{W}$.
Theorem~\ref{t.lusin} then gives some $\tilde{T} \in \mathcal{W} \cap \Homeo(X,\mu)$
such that $d(\tilde{T},T) < \eps$.
This concludes the proof.

\begin{rem}
Notice we haven't used the full strength of theorem \ref{t.meas},
but the mere fact that if $A_*\mu$ is rich then there exists
some $T \in \Aut(X,\mu)$ for which $\LE(A,T)$ is small.
\end{rem}

\subsection{Proof of theorem~\ref{t.contdich}}\label{ss.contdich}

It is interesting to mention that
the following (apparently) stronger form of theorem~\ref{t.cont} holds:

\begin{theorem}\label{t.cont pseudo}
Given $T\in\Homeo(X,\mu)$ and $\eps>0$, there exists $\delta>0$ with the following properties:
Let $A \in C^1(X, \G)$ and assume there exists a
periodic $\delta$-pseudo-orbit $(x_0, \ldots, x_{n-1}, x_n = x_0)$
for $T$ such that:
\begin{itemize}
\item $A(x_{n-1}) \cdots A(x_0)$ is elliptic;
\item $A$ is not locally constant at at least one of the points $x_i$.
\end{itemize}
Then there exists
$\tilde{T} \in \Homeo(X, \mu)$ $\eps$-$C^0$-close to $T$
such that $\LE(A,\tilde{T})$ is arbitrarily close to zero.
\end{theorem}

Theorem~\ref{t.cont pseudo} is actually a corollary of theorem~\ref{t.cont}.
Indeed, given a periodic $\delta$-pseudo-orbit $(x_0, \ldots, x_{n-1}, x_n = x_0)$ for $T$,
there exists a perturbation $\tilde{T}$ of $T$ such that
$\tilde T(x_i) = x_{i+1}$ for $i=0,\ldots,n-1$.
(This follows from \cite[theorem~2.4]{AlPr}, for instance.)

We are going to use the following result due to Avila.
Recall $R_\theta$ denotes the rotation of angle $\theta$,
and $\rho(\mathord{\cdot})$ denotes spectral radius.

\begin{lemma}[Lemma 2 from \cite{Yo}]\label{l.artur}
For every $n\in \N$ and $A_0$,\ldots,$A_{n-1}\in \G$, there is
$\theta \in \R$ such that
$$
R_\theta A_{n-1} \cdots R_\theta A_1 R_\theta A_0
\text{ is elliptic} \quad \text{and} \quad
|\theta| \le \frac{C}{n} \log \rho(A_{n-1} \cdots A_0),
$$
where $C>0$ is some constant.
\end{lemma}

\begin{proof}[Proof of theorem~\ref{t.contdich}]
By proposition~\ref{p.generic fubini}, it suffices to prove
that the set of $(A,T)$ such that either
$(A,T)$ is uniformly hyperbolic or $\LE(A,T)=0$ is generic in
$C^r(X,\G) \times \Homeo(X,\mu)$.
The uniformly hyperbolic cocycles $(A,T)$ form an open set.
Also, the function
$$
\LE : C^r(X,\G) \times \Homeo(X,\mu) \to \R
$$
is an upper semicontinuous function.
So to prove the theorem, we have to show that
if $(A,T)$ is not uniformly hyperbolic then
for every $\eps>0$ there exist
$\tilde{A}$ $C^r$-close to $A$ and $\tilde{T}$ $C^0$-close to $T$
such that $\LE(A,T)<\eps$.

Fix $\eps>0$ and a cocycle $(A,T)$ which is not uniformly hyperbolic.
Making a $C^r$-perturbation if necessary, 
we can assume that $A$ is nowhere locally constant.

Let $\delta>0$ be given by theorem~\ref{t.cont pseudo}.
Because $T$ preserves volume, there is some $\ell_0 \in \N$ such that
for every pair of points $y$, $x \in X$ there exists a
$\delta$-pseudo-orbit $(y_0, y_1, \ldots, y_\ell)$ for $T$,
such that $\ell < \ell_0$, $y_0 = y$, and $y_\ell = x$.

Since $(A,T)$ isn't uniformly hyperbolic,
there exist arbitrarily large $n\in \N$ and $x \in M$ such that
$\| A^n_T(x) \| < (1+\delta)^n$.
Fix $n$ and $x$ with $n > \ell_0/\delta$.

By concatenation we obtain a $\delta$-pseudo-orbit
$(x_0, x_1, \ldots, x_{n+\ell})$
of length $n + \ell < n + \ell_0$ with $x_0 = x_{n+\ell} = x$
and such that
$$
\| A(x_{n+\ell-1}) \cdots A(x_0) \| < \|A\|_\infty^\ell (1+\delta)^n.
$$
According to lemma~\ref{l.artur}, there exists $\theta$
with $|\theta|<\mathrm{const.}\delta$ such that
$$
A(x_{n+\ell-1}) R_\theta \cdots A(x_1) R_\theta A(x_0) R_\theta
$$
is an elliptic matrix.
So theorem~\ref{t.cont pseudo} applies to 
the $C^r$-perturbation $\tilde A = A R_{\theta}$ of $A$,
showing that there exist $\tilde T$ $\delta$-$C^0$-close to $T$ such that
$\LE(\tilde A, \tilde T)$ is as small as we want.
\end{proof}

\section{The discrete case. Questions}\label{s.discrete}

In this section we prove theorem~\ref{t.discrete} and,
in \S~\ref{ss.openquestions}, we pose related problems.

\subsection{Preparation}

\paragraph{Uniformly hyperbolic sets, elliptic products.}

The set of $\Sigma \in \G^N$ which are uniformly hyperbolic is open in~$\G^N$, see \cite{Yo}.
In fact, we will use the following result,
which is indeed an immediate corollary of Avila's lemma~\ref{l.artur}:
\begin{othertheorem}[Proposition 6 in \cite{Yo}]\label{t.hyp or ell}
There is an open and dense subset $\cR_0 \subset \G^N$ such that
if $\Sigma \in \cR_0$ then
either $\Sigma$ is uniformly hyperbolic
or there is an elliptic matrix in the semigroup
$\langle \Sigma \rangle$ generated by~$\Sigma$.
\end{othertheorem}

\paragraph{Liouville pairs.}

Recall that $\rho$ denotes the spectral radius.

\begin{defn}\label{d.liouville}
Let $\psi \colon \N \to \N$, with $\lim_{n \to \infty} \psi(n)=\infty$.
If $R$ and $H$ belong to $\G$, we say that the pair
$(R,H)$ is $\psi$-Liouville if $R$ is elliptic and
$$
\liminf_{n \to +\infty}
\frac{1}{\psi(n)} \log \rho\big(R^n H^{\psi(n)}\big) = 0.
$$
\end{defn}

Notice that if $H$ is not hyperbolic then, for any $\psi$,
$(R,H)$ is $\psi$-Liouville for every elliptic $R$.

\begin{lemma}\label{l.lio is res}
Given any $\psi\colon  \N \to \N$  with $\lim_{n \to \infty} \psi(n)=\infty$,
let $\cR$ be the set of $(R,H)\in \G^2$
such that $R$ is not elliptic or
$(R,H)$ is $\psi$-Liouville.
Then $\cR$ is a residual subset of $\G^2$.
\end{lemma}

\begin{proof}
Let $\G_\mathrm{ell}$ be the subset of $\G$ formed by elliptic matrices,
and let $\cL \subset \G^2$ be the set of $\psi$-Liouville pairs.
We have $\cL = \bigcap_{m\ge 1, \eps>0} U_{m,\eps}$, where
$$
U_{m,\eps}=
\left\{(R,H)\in \G_\mathrm{ell} \times \G; \;
\exists n \ge m \text{ s.t.~}
\frac{1}{\psi(n)} \log \rho\big(R^n H^{\psi(n)}\big) <\eps \right\}.
$$
Each $U_{m,\eps}$ is open and we have to show it is dense in $\G_\mathrm{ell} \times \G$.
Given $(R,H)\in \G_\mathrm{ell} \times \G$, with $H$ hyperbolic,
take a basis of $\R^2$ such that we can write
$$
H=
\begin{pmatrix}
\lambda & 0 \\ 0 & \lambda^{-1}
\end{pmatrix},
\quad |\lambda|>1.
$$
Arbitrarily close to $R$, there is an elliptic matrix $\tilde{R}$ such that
$\tilde{R}^n(1,0) \in \R(0,1)$ for some $n \ge m$, that we can choose satisfying $1< e^{\eps \psi(n)}$.
Hence
$$
\tilde{R}^n=
\begin{pmatrix}
0 & c \\ b & d
\end{pmatrix}
\quad \text{and} \quad
\tilde{R}^n H^{\psi(n)}=
\begin{pmatrix}
0 & c \lambda^{-\psi(n)} \\ b \lambda^{\psi(n)} & d \lambda^{-\psi(n)}.
\end{pmatrix}
$$
If $n$ is chosen large enough we have
$$
\big\lvert \trace \tilde{R}^n H^{\psi(n)} \big\rvert =
|d| |\lambda|^{-\psi(n)} < \mathrm{Const} \cdot \|\tilde R\| \; |\lambda|^{-\psi(n)} < 2.
$$
Therefore $\rho\big(\tilde{R}^n H^{\psi(n)}\big)=1<e^{\eps \psi(n)}$,
that is, $(\tilde{R},H) \in U_{m,\eps}$.
\end{proof}

\paragraph{Monomials.}

Let ${\mathcal N} = \{1,2,\ldots,N\}$.
To every word $(k_1,\ldots,k_\ell) \in {\mathcal N}^\ell$, $\ell\geq 1$,
we can associate a map 
$F\colon \G^N \to \G, (A_1,\ldots,A_N)\mapsto A_{k_1}\cdots A_{k_\ell}$,
which is called a \emph{monomial}.
For each $i \in {\mathcal N}$, let us write
$$
m_i(F) = \# \big\{ j\in \{1,\ldots,\ell\}; \; k_j=i \big\},
$$
that is, the number of appearances of the letter $A_i$ in the monomial $F$.
Let us call two monomials $F_1, F_2\colon  \G^N \to \G$ \emph{independent}
if the vectors $(m_1(F_1), \ldots, m_N(F_1))$ and $(m_1(F_2), \ldots, m_N(F_2)) \in \R^N$
are non-collinear.

\begin{lemma} \label{l.monomials}
Let $F_1, F_2 \colon \G^N \to \G$ be independent monomials,
and let $F = (F_1, F_2)\colon  \G^N \to \G^2$.
Then for every residual subset $R$ of $\G^2$, the
set $F^{-1}(R)$ is residual in $\G^N$.
\end{lemma}

\begin{proof}
Let $C \subset \G^N$ be the set of critical points of $F$.
We will show that $C$ has empty interior.
That implies the lemma, because
$F$ restricted to the open dense set $\G^N \setminus C$ is an open map.

The derivative of $F$ at $(\id,\ldots,\id)$ is easily computed; it is:
$$
(a_1, \ldots, a_N) \in \mathit{sl}(2,\mathbb{R})^N \mapsto
\left( \sum_{i=1}^N m_i(F_1) a_i, \sum_{i=1}^N m_i(F_2) a_i \right)
\in \mathit{sl}(2,\mathbb{R})^2.
$$
Due to the independence assumption, $DF(\id,\ldots,\id)$ is surjective,
that is, $(\id,\ldots,\id)\notin C$.
Assume $C$ has an interior point $x$.
Consider a real-analytic path $[0,1] \to \G^N$ from $x$ to $(\id,\ldots,\id)$.
Bearing in mind that $C$ is the zero-set of some real-analytic function,
we reach a contradiction.
\end{proof}

\subsection{Proof of theorem~\ref{t.discrete}}
In all the proof we fix the function $\psi(n)=n$.

First we define the residual set $\cR \subset \G^N$ for which we will prove the conclusion of the theorem.
Given two independent monomials $F_1$, $F_2\colon  \G^N \to \G$,
let $\cR(F_1,F_2)$ be the set of all $\Sigma \in \G^N$ such that
$$
\text{
$F_1(\Sigma)$ is not elliptic or $(F_1(\Sigma),F_2(\Sigma))$ is $\psi$-Liouville.
}
$$
By lemmas~\ref{l.lio is res} and \ref{l.monomials},
$\cR(F_1,F_2)$ is a residual subset of $\G^N$.
Take the intersection over all independent pairs $F_1$, $F_2$
and call it $\cR_1$.
Finally, let $\cR = \cR_0 \cap \cR_1$,
where $\cR_0$ is the set from theorem~\ref{t.hyp or ell}.

From now on fix some $\Sigma \in \cR$.
If $\Sigma$ is uniformly hyperbolic, there is nothing to do.
In the other case, since $\Sigma \in \cR_0$,
there is a monomial $F_1$ such that $R = F_1(\Sigma)$ is elliptic.
$F_1$ will be fixed from now on.
By construction, $(F_1(\Sigma), F_2(\Sigma))$ is $\psi$-Liouville for every monomial $F_2$
which is independent from~$F_1$.

Let $A\colon  X \to\Sigma$ be a measurable function such that every matrix in $\Sigma$
is attained on a positive measure set of $X$.
As usual, we assume $X$ is the unit interval $\I$.
We can also suppose there is a partition $\I = I_1 \sqcup \cdots \sqcup I_N$ into
intervals such that $A|_{I_i} = A_i$,
where $\Sigma = (A_1,\ldots,A_N)$.

Now let $T\colon \I \to \I$ be any given automorphism.
We will explain how to perturb $T$ in the weak topology to make the exponent small.
By proposition~\ref{p.semiT}, this will conclude the proof.

\smallskip

Using theorem~\ref{t.halmos}, we may begin with $T$ equal to a cyclic interval permutation of
some arbitrarily high rank $M$.

We will of course perturb $T$ further, but will work only with automorphisms
that are (not necessarily cyclic) interval permutations.
In this regard, a sequence of \emph{disjoint} intervals $J_i = T^i(J_1), i=0,\ldots,\ell-1$,
is called a \emph{tower of height $\ell$}.
The tower is said to be \emph{cyclic} if in addition $T^\ell(J_1) = J_1$.
If moreover the map $A\colon \I \to \Sigma$ is constant on each interval $J_i$
then we can talk about the \emph{product of matrices along the tower},
that we denote by $A(J_M) \cdots A(J_1)$.

\smallskip

Since the rank $M$ is high, most of the intervals
$\left[ \frac{j}{M}, \frac{j+1}{M} \right)$, $j=0,1,\ldots,M-1$,
will be completely contained in one
of the intervals $I_i$ (where $A$ is constant).
By changing $T$ on a set of small measure,
we may assume the collection of those ``good'' intervals is cyclically permuted by $T$.
The union of the ``bad'' intervals is now an invariant set of small (less than $N/M$) measure,
and so its contribution to the mean Lyapunov exponent is small.
So, to simplify writing, we will assume that all intervals are good.

\smallskip

Write $F_1 (A_1,\ldots,A_N) = A_{k_p} \cdots A_{k_1}$.
Among the intervals that are permuted by $T$, select some
$J_1$, \ldots, $J_p$ such that $A|_{J_i} = A_{k_i}$.
Since $M$ can be chosen arbitrarily high compared to $p$,
the measure of $\cT_1 = \bigsqcup_{i=1}^p J_i$ is small.
So, after another perturbation, we can assume that
that the dynamics of $T$ decomposes into two cyclic towers,
the smaller of which is
$$
J_1 \to J_2 \to \cdots \to J_p \to J_1.
$$
Call this tower $\cT_1$.
The product of matrices along it is precisely the elliptic matrix $F_1(\Sigma)$.

Let $\cT_2$ be the other, bigger, tower.
Consider the product of matrices along $\cT_2$;
as a function of $\Sigma$ that product is by definition a monomial $F_2(\Sigma)$.
We may assume that $F_1$ and $F_2$ are independent monomials.
In fact, via a small perturbation of $T$
we can remove a single level of the bigger tower
to make $F_1$ and $F_2$ independent.
The removed interval becomes an invariant set of small measure and can be disregarded.

Hence, by definition of $\cR$, the pair
$(R,H) = (F_1(\Sigma), F_2(\Sigma))$ is $\psi$-Liouville.
That is, there is an integer $n$ such that
\begin{equation} \label{e.psi}
\frac{1}{n}\log \rho(H^n R^n) < \eps,
\end{equation}
where $\eps>0$ is any fixed small number.

Decompose each interval $J_i$ into $n$ intervals of equal length
$J_i = J_{i,1} \sqcup J_{i,2} \sqcup \cdots \sqcup J_{i,p}$.
Modify slightly $T$ in order to form the following tower of height $np$:
\begin{align*}
&J_{1,1} \to \cdots \to J_{p,1} \to \\
&J_{1,2} \to \cdots \to J_{p,2} \to \\
&\cdots \to \\
&J_{1,n} \to \cdots \to J_{p,n} \to J_{1,1}.
\end{align*}
The product along this new tower is $R^n$.
In the same way we decompose the $\cT_2$ tower in $n$ towers
that we unfold as above into a single tower along which the product of matrices will be $H^n$.
As sets, the two new towers are still the same $\cT_1$ and $\cT_2$.

By our construction, $\cT_1$ and $\cT_2$ have bases of the same size.
So we can actually concatenate them one on top of the other to get
a single cyclic tower along which the matrix product is $H^n R^n$.
(This is done by composing on the left the dynamics with a map that permutes
the bases of the towers.)
Since almost all the space is covered by this tower, we deduce from~\eqref{e.psi}
that the integrated Lyapunov exponent corresponding to the perturbed dynamics is small.
This proves theorem~\ref{t.discrete}.

\subsection{Some open questions}\label{ss.openquestions}

\begin{problem}\label{pb.exist}
Does a finite set $\Sigma \subset \G$ with the following properties exist?
\begin{enumerate}
\item $\Sigma$ cannot be approximated by a uniformly hyperbolic set;
\item there exists a measurable map $A \colon \I \to \Sigma$ which assumes
every value in $\Sigma$ on a set of positive measure
such that $\LE(A,T) > 0$ for \emph{every} $T \in \Aut(\I,m)$.
(Or even stronger, such that $\LE(A,T)\ge\lambda_0>0$ for every $T$.)
\end{enumerate}
\end{problem}

By theorem~\ref{t.discrete}, those $\Sigma$ with $\#\Sigma = N$
form a meager subset of $\G^N$.

A positive answer to the following more elementary question would,
by the ergodic theorem, answer problem~\ref{pb.exist} (in the stronger form) positively:
\begin{problem}\label{pr.elementary}
Does there exist a pair of matrices $A_1$, $A_2\in \G$ with the following properties?
$A_1$ is hyperbolic, $A_2$ is elliptic,
and there are constants $0<p<1$ and $\lambda>1$ such that for every word
$A_{i_1} A_{i_2} \cdots A_{i_n}$ satisfying the \emph{frequency condition}
\begin{equation}\label{e.frequency}
\# \big\{j\in\{1,\ldots,n\}; \; i_j=2 \big\} < pn \, ,
\end{equation}
we have
\begin{equation}\label{e.growth}
\| A_{i_1} A_{i_2} \cdots A_{i_n} \| > \lambda^n \, .
\end{equation}
\end{problem}

Fixing some integer $N\ge 2$, we can also ask
whether the set of $\Sigma \in \G^N$  that have the properties as in problem~\ref{pb.exist}
has positive, or even full measure in the complement of the hyperbolicity locus in $\G^N$.

\begin{rem}
From lemma~\ref{l.lio is res} we see that even if the right hand side
in (\ref{e.frequency}) is replaced by any function $\phi(n)$
such that $\phi(n) \to \infty$ then the set of  $\Sigma \in \G^2$ that satisfy the conclusion of
problem~\ref{pr.elementary} is meager.
\end{rem}

\appendix
\section{Appendices}

\subsection{Derived cocycles}

Given a set $W \subset X$ of positive measure,
we define the \emph{first return map} $T_W \colon W \to W$
by $T_W(x) = T^{n_W(x)}(x)$, where $n_W(x) = \min \{n \ge 1; \; T^n(x) \in W \}$.
$n_W$ and $T_W$ are defined a.e.~and we have $T_W \in \Aut(W, \mu|_W)$.

Define $A_{T,W} \colon  W \to \G$ as $A_{T,W}(x) = A_T^{n_W(x)}(x)$.
The pair $(T_W, A_{T,W})$ is a $\G$-cocycle on $(W, \mu|_W)$.
This is called the \emph{derived} or \emph{induced cocycle}.

\begin{prop} \label{p.derived}
We have $\LE(A_{T,W}, T_W) = \LE(A|_{\cT}, T|_{\cT})$, where
$\cT$ is the $T$-invariant set $\bigcup_{n\in \Z} T^{-n}(W)$.
\end{prop}

In fact, if $T$ is ergodic then $\cT = X \bmod{0}$ and the proposition is lemma~2.2 from~\cite{K}.
It is easy to adapt that proof to the general case, using Kac's formula in the form:
$$
\int_W n_W \; d \mu = \mu(\cT).
$$

\subsection{Semicontinuity}

It is well-known that:
$$
\LE(A,T) = \inf_N \frac{1}{N} \int_X \log\|A_T^N\| \, d\mu.
$$
Among the consequences, we have semicontinuity of $\LE (A, \cdot)$:

\begin{prop}\label{p.semiT}
Let $A\colon X \to \G$ be measurable and such that $\log \|A\| \in L^1(\mu)$.
Then the function $T \in \Aut(X,\mu) \mapsto \LE(A,T) \in \R$ is upper semicontinuous.
\end{prop}

\begin{proof}
We may assume that $X$ is the unit interval $\I$ and $\mu$ is Lebesgue measure $m$.
The weak topology in $\Aut(\I,m)$ is then given by the \emph{weak metric}
$$
d (S,T) = \inf \{\rho>0 ; \; m(\{|S-T|>\rho\}) < \rho \}.
$$

Let $A\colon \I\to\G$, $T\in \Aut(\I,m)$, and $\eps>0$ be fixed.
There exists $N\in\N$ such that
$$
\LE(A,T) > - \eps + \frac{1}{N} \int_\I \log\|A_T^N\| \, dm.
$$
Since $\log\|A\|$ is integrable, there is $\delta_1>0$ such that if $Z \subset \I$
has measure $m(Z)<\delta_1$ then $\int_Z \log\|A\|\, dm< \eps$.
By Lusin's theorem, there exists a compact set $K \subset \I$ such that
the functions $A|_K$ and $T|_K$ are continuous, and $m(K^c) < \delta_1/(2N)$.
Let $C = \sup_K \|A\|$.
There is $\delta_2>0$ such that if
$A_1,\ldots,A_N, B_1, \ldots, B_N \in \G$ are matrices with norm at most $C$ and
$\|A_i - B_i\|< \delta_2$ for each $i$ then $\|\prod_N^1 A_i - \prod_N^1 B_i\| < \eps$.
Let $\delta_3>0$ be such that if $x,y \in K$, $|x-y| < \delta_3$ then
$\|A(x)-A(y)\|<\delta_2$.
Take numbers $\eta_1 > \cdots > \eta_{N-1} >0$ such that $\eta_1 = \delta_3/2$ and
$$
x,y \in K, \ |x-y|< 2\eta_{i+1} \ \Rightarrow \
|T(x) - T(y)| < \eta_i.
$$
Let $\rho = \min(\eta_{N-1}, \delta_1/(2N))$.

Now assume $S\in\Aut(\I,m)$ is such that $d(S,T)<\rho$.
Let $W = \{|S-T| \le \rho\}$; then $m(W^c) < \rho$.
Define
$$
G = \bigcap_{i=0}^{N-1} \left[ T^{-i} (K \cap W)\cap S^{-i} (K \cap W) \right].
$$
Then $G^c$ has small measure: $m(G^c) \le N m(K^c + W^c) < \delta_1$.
We are going to bound the expression
$\frac{1}{N}\int_\I \log \| A_S^N\| \, dm$.
To do so, we are going to split the integral in two parts,
$\int_\I = \int_{G^c} + \int_G$.
For the first part, we have
$$
\frac{1}{N}\int_{G^c} \log \| A_S^N \| \, dm
\le \frac{1}{N}\sum_{i=0}^{N-1} \int_{S^i(G^c)} \log \| A \|  \, dm
\le \eps.
$$
For the second part,
\begin{align*}
\frac{1}{N}\int_G \log \| A_S^N \| \, dm
         &\le \frac{1}{N} \int_G \log\|A_T^N\| \, dm + \frac{1}{N} \int_G \|A_S^N - A_T^N\| \, dm \\
         &<   \LE(A,T) + \eps + \frac{1}{N} \int_G \|A_S^N - A_T^N\| \, dm \, .
\end{align*}
(We used that $\log(a+b) \le \log a + b$ for $a\ge 1$, $b\ge 0$.)

Let $x\in G$.
We claim that $|T^i(x)-S^i(x)| \le 2 \eta_{N-i}$ for all $i=0,1,\ldots,N-1$.
This is easily shown by induction:
\begin{align*}
|T^{i+1}(x)- S^{i+1}(x)| &\le |T(T^i(x))- T(S^i(x))| + |T(S^i(x)) - S(S^i(x))| \\
                         &\le  \eta_{N-i-1} + \rho
                         \le  2\eta_{N-i-1}.
\end{align*}
In particular, for all $i$ we have $|T^i(x)-S^i(x)| < \delta_3$
and thus $\|A(T^i(x)) - A(S^i(x))\| < \delta_2$.
Therefore $\|A_T^N(x) - A_S^N(x) \| < \eps$.

Summing the two parts, we conclude that
$$
\LE(A,S) \le \frac{1}{N}\int_G \log \| A_S^N \| \, dm
         \le \LE(A,T) + \eps + \frac{\eps}{N} \, .
$$
This shows upper semicontinuity.
\end{proof}

\begin{rem}
For the semicontinuity of $\LE(\cdot,T)$ in the $L^1$-topology, see~\cite{ArB}.
\end{rem}

\subsection{An auxiliary result from measure theory}

The aim of this section is to establish the lemma~\ref{l.thouvenot} below,
which we use a few times in section~\ref{s.proofmain}.
In what follows, $\I$ denotes the unit interval $[0,1]$,
$m$ denotes Lebesgue measure in $\I$ or (by abuse of notation) in the square $\I^2$,
and $\pi\colon \I^2\to\I$ the projection in the first coordinate.

\begin{lemma}\label{l.thouvenot}
Let $Y$ be either $\G$ or $\P^1$.
Let $A, A'\colon  \I \to Y$ be measurable functions such that
$A_* m = A'_* m = \nu$.
Then there exists $S\in \Aut(\I^2,m)$ such that $A' \circ \pi \circ S = A \circ \pi$ $m$-a.e.
\end{lemma}

The lemma is a straightforward consequence of the work of Rokhlin~\cite{R},
as we explain below.

Let $(X,\mu)$ and $(Y,\nu)$ be Lebesgue spaces and let
$h\colon  (X,\mu) \to (Y,\nu)$ be a homomorphism (that is, a measurable map with $\nu = h_* \mu$).
Consider the Rokhlin disintegration of the measure $\nu$ along fibers of $h$:
For $\nu$-a.e.~$y \in Y$ we have a probability measure $\mu_y$ on the set $h^{-1}(y)$.

Lemma~\ref{l.thouvenot} is in fact a particular case of the following:

\begin{prop}\label{p.generalthouvenot}
Let $h, h'\colon  (X,\mu) \to (Y,\nu)$ be two homomorphisms between Lebesgue spaces.
Let $\{\mu_g\}$ and $\{\mu'_g\}$ be the respective Rokhlin disintegrations of the measure $\mu$.
Assume that for $\nu$-a.e.~$y\in Y$,
the measures $\mu_y$ and $\mu'_y$ are non-atomic.
Then there exists an automorphism $S\colon (X,\mu) \to (X,\mu)$ such that $h' \circ S = h$ $\mu$-a.e.
\end{prop}

\begin{proof}
Consider the two decompositions
$\zeta = \{h^{-1}(y)\}_{y\in Y}$ and $\zeta' = \{(h')^{-1}(y)\}_{y\in Y}$
of the space $X$.
The factor spaces $X/\zeta$ and $X/\zeta'$ are isomorphic
(to $(Y,\nu)$, see \cite[p.~32]{R}).
Therefore the result of \cite[p.~51]{R},
together with the assumption that $\mu_y$ and $\mu'_y$ are non-atomic,
gives that the decompositions $\zeta$ and $\zeta'$ are isomorphic mod~$0$.
This means (see \cite[p.~9]{R}) that there exist isomorphisms $U$ and $V$
which make the diagram below commutes mod $0$:
$$
\xymatrix{X       \ar[r]^U \ar[d] & X        \ar[d] \\
          X/\zeta \ar[r]_V        & X/\zeta' }
$$
(Vertical arrows denote quotient maps.)
From $U$ we construct the desired automorphism $S$.
\end{proof}

\subsection{Uniform spectral radius theorem}\label{ss.unif spectral}

Let $\|\mathord{\cdot}\|$ be an operator norm in $\mathit{GL}(d,\R)$,
and let $\rho(A)$ denote the spectral radius of $A\in \mathit{GL}(d,\R)$.

\begin{prop}\label{p.unif spectral}
We have $\rho(A) = \lim_{n\to+\infty} \|A^n\|^{1/n}$ for every $A \in \mathit{GL}(d,\R)$,
and the convergence is uniform in compact subsets of $\mathit{GL}(d,\R)$.
\end{prop}

\begin{proof}
The first part is the spectral radius theorem.
Now, fix $A_0 \in \mathit{GL}(d,\R)$ and $\eps>0$.
Let $n_0 \in \N$ be such that $\|A_0^{n_0}\|^{1/n_0} < \rho(A_0)+\eps$.
Let $\delta>0$ be such that $\|A-A_0\|<\delta$ implies
$\|A^{n_0}\|^{1/n_0} < \rho(A_0)+ 2\eps$ and $\rho(A) > \rho(A_0)-\eps$.
For $n \in \N$, we write $n=mn_0+k$ with $m\ge 0$ and $0\le k < n_0$.
If $n$ is sufficiently large and $A$ is $\delta$-close to $A_0$ then
\begin{multline*}
\rho(A) \le \|A^n\|^{1/n} 
\le \|A^{n_0}\|^{m/n} \; \|A\|^{k/n} \\
< (\rho(A_0)+2\eps)^{mn_0/n} (\|A_0\|+\delta)^{k/n}
< \rho(A_0)+3\eps
< \rho(A)+4\eps.
\end{multline*}
This shows uniform convergence in the ball $\{A; \; \|A-A_0\|<\delta\}$.
\end{proof}

\subsection{Generic measurable dichotomy for $L^\infty$ or $C^0$ cocycles}\label{ss.easydich}

Here we show proposition~\ref{p.easy dich}.
For that we need the proposition below, which is also used
in the proof of theorem~\ref{t.contdich}.

\begin{prop}\label{p.generic fubini}
If $X$, $Y$ are separable Baire spaces and $\cR \subset X \times Y$ is residual,
then there is a residual subset $\cR'\subset X$ such that for every $x \in \cR'$,
the set $\cR_x = \{y \in Y;\; (x,y)\in\cR\}$ is a residual subset of $Y$.
\end{prop}

\begin{proof}
First let $A \subset X \times Y$ be an open and dense set.
For any open set $V \subset Y$, let
$X_{A,V} = \{x \in X; \; \text{there exists $y \in V$ such that $(x,y) \in A$}\}$.
Then $X_{A,V}$ is open and dense in $X$.
Let $\mathcal{V}$ be a countable base of open subsets of $Y$,
and consider the residual set $X_A = \bigcap_{V \in \mathcal{V}} X_{A,V}$.
If we define $A_x = \{y \in Y; \; (x,y) \in A\}$,
then $A_x$ is open and dense in $Y$ for every $x \in X_A$.

Now, given a residual set $\cR \subset X \times Y$,
write $\cR = \bigcap_{n \in \N} A_n$, where $A_n$ are open and dense.
Let $\cR' = \bigcap_{n \in \N} X_{A_n}$.
Then for every $x\in \cR'$, the fiber
$\cR_x = \bigcap_{n \in \N} (A_n)_x$ is a residual subset of $Y$.
\end{proof}

\begin{proof}[Proof of proposition~\ref{p.easy dich}]
We will only prove the $L^\infty$ statement, because the $C^0$ one is analogous.
Note that if $A\in L^\infty(X,\G)$ then one should read the first alternative
in the measurable dichotomy with ``a.e.~points'' in place of ``all points''.

By propositions~\ref{p.semiT} and~\ref{p.generic fubini},
we only have to show that if the essential image of $A\in L^\infty(X,\G)$ is not
an uniformly hyperbolic set (see definition~\ref{d.unif hyp}) then for every
$T\in \Aut(X,\mu)$ there exist $\tilde{A}$ and $\tilde{T}$ close to $A$ and $T$ respectively in
the $L^\infty$ and weak topologies, such that $\LE(\tilde{A},\tilde{T})$ is small.

Fix such an $A$.
Given $T\in \Aut(X,\mu)$, 
to say that the cocycle $(A,T)$ is \emph{not} uniformly hyperbolic means that
$$
\forall c>0, \ \forall \lambda>1 , \ \exists n \in \N \text { s.t. }
\mu \left\{ x; \; \|A_T^n(x)\| < c \lambda^n \right\} >0 \, .
$$
Therefore the set of $T\in \Aut(X,\mu)$ such that $(A,T)$ is not uniformly hyperbolic
is a $G_\delta$ set (in the weak topology).
That set is also dense, by Baire and 
the assumption that the essential image of $A$ is not an uniformly hyperbolic set.

So given any automorphism $T$, we can find a weak perturbation 
$\tilde{T}$ which is ergodic (see remark~\ref{r.generic ergodicity}) and such that
$(A,\tilde{T})$ is not uniformly hyperbolic.
By a result from \cite{B} (mentioned in remark~\ref{r.when A varies}),
there exists $\tilde{A}$ $L^\infty$-close to $A$ such that
$\LE(\tilde{A},\tilde{T})=0$.
\end{proof}

\paragraph{Acknowledgements:}
J.B.~thanks IMPA, and universities of Villetaneuse, Orsay, and Jussieu
for their hospitality.
Discussions with
Fr\'ed\'eric Le Roux,           
Ivan Pan,                       
Jorge Vit\'orio Pereira,        
and specially Artur Avila       
and Jean-Paul Thouvenot         
contributed to this paper.
We also thank the referee for his/her suggestions.


\vfill

\noindent
\begin{minipage}{5cm}
{\footnotesize Jairo Bochi\\
{\tt jairo{\@@}mat.ufrgs.br}\\
Instituto de Matem\'atica -- UFRGS\\
Porto Alegre, Brazil
}
\end{minipage}
\hfill
\begin{minipage}{5cm}
{\footnotesize Bassam R.~Fayad \\
{\tt fayadb{\@@}math.univ-paris13.fr} \\
Paris 13\\
Villetaneuse, France
}
\end{minipage}

\end{document}